\documentclass{amsart}
\newcommand{\RN}[1]{\textup{\uppercase\expandafter{\romannumeral#1}}}
\usepackage{amsmath}
\usepackage{amssymb}
\usepackage{amsthm}
\usepackage{enumerate}
\usepackage[hidelinks]{hyperref}
\usepackage[capitalise]{cleveref}
\usepackage{mathtools}
\usepackage{thmtools} 
\numberwithin{equation}{section}

\newcommand\mc{\mathcal}
\newcommand\mf{\mathfrak}
\newcommand\mb{\mathbb}
\newcommand{\rank}{\mathrm{rank}}
\newcommand{\Bl}{\mathrm{Bl}}
\newcommand{\pr}{\mathrm{pr}}
\newcommand{\taut}{\mathrm{taut}}
\newcommand{\D}{\mathrm{D}}
\crefname{equation}{}{}
\newtheorem{theorem}{Theorem}[section]
\newtheorem{lemma}[theorem]{Lemma}

\newtheorem{proposition}[theorem]{Proposition}

\title{Moduli Space of Sheaves and Categorified Commutator of Functors}
\author{Yu Zhao}
\address{Beijing Institute of Technology, Haidian, Beijing, China}
\email{zy199402@live.com}
\theoremstyle{definition}

\newtheorem{definition}[theorem]{Definition}
\newtheorem{example}[theorem]{Example}

\theoremstyle{remark}
\newtheorem{remark}[theorem]{Remark}
\usepackage{epigraph}
\usepackage{comment}
\usepackage{tikz}
\usepackage{tikz-cd}

\usepackage{enumitem}

\begin{document}
\begin{abstract}
  Negu\c{t} constructed an action of the quantum toroidal algebra on the K-theory of the smooth moduli space of stable sheaves over an algebraic surface, generalizing the action studied by Nakajima, Grojnowski, and Baranovsky in cohomology. In this paper, we construct a weak categorification of Negu\c{t}'s action. The new ingredients are two intersection-theoretic descriptions of the quadruple moduli space of stable sheaves.
\end{abstract}
\maketitle

\section{Introduction}

\subsection{Motivation}  Let $S$ be a smooth projective algebraic surface over $\mb{C}$. Let $S^{[n]}$ be the Hilbert scheme of $n$ points on $S$. A novel construction of Nakajima \cite{nakajima1997} and Grojnowski \cite{grojnowski1996} is an action
\begin{equation}
\label{eqn:heis intro}
\text{Heisenberg algebra} \curvearrowright \bigoplus_{n=0}^{\infty} H^*(S^{[n]})
\end{equation}
which realizes the cohomology groups of Hilbert schemes of points as the Fock space representation of the (super)Heisenberg algebra. 

This construction was generalized by Negu\c{t} \cite{negut2017shuffle,neguct2018hecke} in two directions: first, the Hilbert scheme of points was replaced by the moduli space of stable sheaves over $S$; second, the cohomology groups were replaced by the $K$-theory groups. More precisely, under the assumptions in \cref{sec:1.2}, Negu\c{t} constructed a quantum toroidal algebra $\mathrm{U}_{q_{1},q_{2}}(\ddot{\mathrm{gl}}_{1})$ action on the Grothendieck group of $\mc{M}$, which also generalized the construction of Schiffmann--Vasserot \cite{schiffmann2013,MR3150250} when $S=\mb{P}^{2}$ and $\mc{M}$ is the moduli space of framed sheaves and that of Feigin--Tsymbaliuk \cite{MR2854154} when $S=\mb{A}^{2}$ and $\mc{M}$ is the Hilbert scheme of points.

The categorification of the quantum toroidal algebra action was also studied in \cite{neguct2018hecke}, where the relations of the $e_{i}$ actions were lifted to derived categories and Negu\c{t} obtained a weak categorification for the positive part of the quantum toroidal algebra. The purpose of this paper is to construct a weak categorification of Negu\c{t}'s action for the whole quantum toroidal algebra.

\subsection{The main results}
  \label{sec:1.2}
 Fix $(r,c_{1})\in \mb{N}^{>0}\times H^{2}(S,\mb{Z})$ and an ample line bundle $H$. Let $\mc{M}$ be the moduli space of Gieseker $H$-stable sheaves $\mc{F}$ such that $\rank(\mc{F})=r$ and $c_{1}(\mc{F})=c_{1}$, and let $\mc{M}_{k}$ be the subspace of $\mc{M}$ consisting of sheaves $\mc{F}$ such that $c_{2}(\mc{F})=k$, with a universal sheaf $\mc{U}_{k}$ over $\mc{M}_{k}\times S$ for every integer $k$.  Let $q:=K_{S}$ be the canonical line bundle of $S$. We make the following assumptions:
\begin{equation*}
  \gcd(r, c_{1}\cdot H) = 1,  \qquad  q \cong \mc{O}_S \text{ or } c_1(q) \cdot H < 0.
\end{equation*}
 
 It is proved in \cite{negut2017shuffle, neguct2018hecke} that if the above assumptions are satisfied, then
\begin{enumerate}[leftmargin=*]
\item  For any non-split short exact sequence
  $$0\to \mc{F}\to \mc{F}' \to \mb{C}_{x}\to 0,$$
  the sheaf $\mc{F}'\in \mc{M}_{k}$ if and only if $\mc{F}\in \mc{M}_{k+1}$.
\item There exists a constant $\mathrm{const}$, depending only on $S,H,r,c_{1}$, such that $\mc{M}_{k}$ is a smooth projective scheme of dimension $\mathrm{const}+2rk$.
\end{enumerate}

Let $\mf{Z}_{k,k+1}$ be the closed subscheme of $\mc{M}_{k}\times \mc{M}_{k+1}\times S$ consisting of pairs
\begin{equation}
  \label{e11}
  \{(\mc{F}_{-1}\subset_{x} \mc{F}_{0})|\mc{F}_{0}\in \mc{M}_{k}, x\in S\},
\end{equation}
where $\mc{F}_{-1}\subset_{x}\mc{F}_{0}$ means $\mc{F}_{-1}\subset \mc{F}_{0}$ and $\mc{F}_{0}/\mc{F}_{-1}\cong \mb{C}_{x}$. There is a tautological line bundle $\mc{L}$ on $\mf{Z}_{k,k+1}$, such that for each closed point $(\mc{F}_{0}\subset_{x} \mc{F}_{1})$, the fiber of $\mc{L}$ at this point is $\mc{F}_{1}/\mc{F}_{0}$. We define the Fourier-Mukai transform associated to a surface in \cref{FMS}. Given an integer $i$, we consider the Fourier-Mukai kernels $e_{i},f_{i}\in D^{b}(\mc{M}\times \mc{M}\times S)$ as the disjoint union of
\begin{equation*}
  \mc{L}^{i}\mc{O}_{\mf{Z}_{k,k+1}}\in D^{b}(\mc{M}_{k}\times \mc{M}_{k+1}\times S) \text{ and }  \mc{L}^{i-r}\mc{O}_{\mf{Z}_{k,k+1}}\in D^{b}(\mc{M}_{k+1}\times \mc{M}_{k}\times S)
\end{equation*}
 for all $k$, respectively. Let $\Delta_{S}:\mc{M}\times S\to \mc{M}\times \mc{M}\times S\times S$ be the diagonal embedding, and let $\iota:\mc{M}\times \mc{M}\times S\times S\to \mc{M}\times \mc{M}\times S\times S$ be the involution morphism that maps $(x,z,s_{1},s_{2})$ to $(x,z,s_{2},s_{1})$.

\begin{theorem}[See \cref{thm61} for the details]
  \label{thm11}
  We have the decomposition
  $$f_{-i}e_{i}\cong R\iota_{*}e_{-i}f_{i}\bigoplus_{a=-r+1}^{0} R\Delta_{S*}(q^{-a}det(\mc{U}_{k})^{-1}\mc{O}_{\mc{M}_{k}\times S})[1-2a-r]$$
and explicit morphisms
  \begin{equation*}
    \begin{cases}
      e_if_j\to R\iota_{*}(f_{j}e_{i}) & i+j>0, \\
      e_if_j\xleftarrow{} R\iota_{*}(f_{j}e_{i}) & i+j<0.
    \end{cases}
  \end{equation*}
whose cones are filtered by combinations of symmetric and exterior powers of the universal sheaf $\mc{U}_{k}$ and its derived dual.
  \end{theorem}

The relations of the $e_{i}$ actions were lifted to derived categories in Theorem 1.1 of \cite{neguct2018hecke}. This paper lifts the commutators of $e_{i}$ and $f_{j}$ actions to derived categories and thus can be regarded as a weak categorification of Negu\c{t}'s action.
The quantum toroidal algebra can also be presented by the elliptic Hall algebra by Schiffmann \cite{MR2886289} and Burban-Schiffmann \cite{MR2922373}. It contains further commutator relations, which we do not currently know how to categorify. See \cite{neguct2018hecke} for further discussion.

\subsection{Strategy of the proof}
The geometry of nested moduli spaces of sheaves was studied by Negu\c{t} \cite{negut2017shuffle,neguct2017w,neguct2017agt,neguct2018hecke}. We adapt this framework and consider the triple and quadruple moduli spaces $\mf{Z}_{-},\mf{Z}_{+}$ and $\mf{Y}$ in \cref{s13}. Unlike the case $r=1$, which was constructed by the author \cite{zhao2020categorical}, when $r>1$ the moduli space $\mf{Z}_{+}$ is neither equidimensional nor Cohen-Macaulay, and thus the dg/derived enhancement must be considered. In this paper, we obtain new descriptions of $\mf{Y}$ in \cref{p36} through blow-ups, which imply the vanishing property we need.

\subsection{Related work} 
\subsubsection{The categorified Hall algebra} The formulation using derived algebraic geometry is inspired by the work of Porta-Sala \cite{porta2019categorification}, Diaconescu-Porta-Sala \cite{diaconescu2020mckay} and Toda \cite{toda2019categorical,toda2020hall,toda2021categorical,toda2021derived} on the categorified Hall algebra. Toda \cite{toda2021derived} proved a recent conjecture of Jiang \cite{jiang2021derived}, who obtained a semiorthogonal decomposition of the derived category of Grassmannians over a coherent sheaf of cohomological dimension $1$. As an application, Koseki \cite{koseki2021categorical} obtained a categorical blow-up formula for Hilbert schemes of points on a smooth algebraic surface.
\subsubsection{The categorification of the $A_{q,t}$ algebra} The quantum toroidal algebra is a subalgebra of the $A_{q,t}$-algebra, which was invented by Carlsson-Mellit \cite{shuffle} to solve the shuffle conjecture. Gonzalez-Hogancamp \cite{GH} obtained a skein-theoretic formulation of the polynomial representation of the $A_{q,t}$ algebra at $q=t^{-1}$ and obtained a categorification of it using the derived trace of the Hecke category. 

\subsubsection{The loop categorification of quantum loop $sl_{2}$}
Shan-Varagnolo-Vasserot \cite{shan2019coherent} constructed an equivalence of graded Abelian categories from a category of representations of the quiver-Hecke algebra of type $A_{1}^{(1)}$ to the category of equivariant perverse coherent sheaves on the nilpotent cone of type A. This gives a representation-theoretic categorification of the preprojective K-theoretic Hall algebra considered by Schiffmann-Vasserot \cite{schiffmann2017cohomological} when the quiver is type $A_{1}$. For affine type quivers, Varagnolo-Vasserot \cite{Varagnolo2021}  proved that the $K$-theoretic Hall algebra of a preprojective algebra is isomorphic to the positive half of a quantum toroidal quantum group.

\subsubsection{Derived projectivizations of complexes}
After we finished the first version of this paper, Qingyuan Jiang sent us his paper \cite{jiang2022derived}. That paper constructed projectivizations of complexes in the setting of derived algebraic geometry and proved a generalized Serre theorem. The construction of Hecke correspondences was also obtained there in the setting of derived algebraic geometry.

\subsection{Organization of the paper}
In \cref{sec2}, we introduce the background and notation used in this paper and review the relevant material from derived algebraic geometry. In \cref{sec:n}, we recall the nested, triple, and quadruple moduli spaces of stable sheaves and prove new descriptions of the quadruple moduli space of stable sheaves via blow-ups. In \cref{sec6}, we prove the main theorem.

\subsection{Acknowledgements}
The author is grateful to Andrei Negu\c{t} for teaching the geometry of the moduli space of sheaves. The author would like to thank Yukinobu Toda, Mikhail Kapranov, Olivier Schiffmann,  Qingyuan Jiang, Francesco Sala, Mauro Porta, and Jeroen Hekking for many useful discussions. The author is supported by World Premier International Research Center Initiative (WPI initiative), MEXT, Japan and JSPS KAKENHI Number 22K13889.

\section{Notation and background}
\label{sec2}
In this section, schemes are denoted by $X,Y,Z$ and locally free sheaves by $V,W$. 
\subsection{Projectivization of two-term complexes}
Given a cosection of a locally free sheaf $V$ on the scheme $X$, i.e. a morphism $f:V\to \mc{O}_X$, the zero locus $\mc{Z}(f)$ is defined as the spectrum of $\mathrm{coker}(f)$, and the Koszul complex $\wedge^{\bullet}f$ is defined as the complex
   $$0\to \wedge^{\rank(V)}V\to \cdots \to V\xrightarrow{f}\mc{O}_{X}\to 0,$$
and we say that $f$ is regular if $\wedge^{\bullet}f$ is a resolution of $\mc{O}_{\mc{Z}(f)}$. 
The following lemma is well-known:
\begin{lemma}
  \label{l21}
If $X$ is smooth and $\dim(\mc{Z}(f))=\dim(X)-\dim(V)$, then $\mc{Z}(f)$ is a complete intersection variety; hence $f$ is regular.  
\end{lemma}

Given a coherent sheaf $U$ on $X$, the projectivization $\mb{P}_{X}(U)$ is a scheme over $X$ with the following functorial property that for any $X$-scheme $f:T\to X$,
$$\mathrm{Hom}_{X}(T,\mb{P}_{X}(U))=\{f^{*}U\twoheadrightarrow E|E\in \mathrm{Pic}(T) \text{ and } f^{*}U\twoheadrightarrow E \text{ is surjective}\}.$$
Let $\pr_{U}:\mb{P}_{X}(U)\to X$ be the projection morphism, and let $\mc{O}_{\mb{P}_{X}(U)}(1)$ be the tautological line bundle on $\mb{P}_{X}(U)$. For a two-term complex $U:=\{W\xrightarrow{s} V\}$ of locally free sheaves over $X$, we abuse notation and write
\begin{equation*}
\mb{P}_{X}(U):=\mb{P}_{X}(\mathrm{coker}(s)) \text{ and } \mc{O}_{\mb{P}_{X}(U)}(1):=\mc{O}_{\mb{P}_{X}(\mathrm{coker}(s))}(1).
\end{equation*}
Let $\taut_{s}:\pr_{V}^{*}W\otimes \mc{O}_{\mb{P}_{X}(V)}(-1)\to \mc{O}_{\mb{P}_{X}(V)}$ be the following composition of morphisms of locally free sheaves on $\mb{P}_{X}(V)$:
 \begin{equation*}
    \pr_{V}^{*}W\otimes \mc{O}_{\mb{P}_{X}(V)}(-1)\xrightarrow{\pr_{V}^{*}s\otimes \mc{O}_{\mb{P}_{X}(V)}(-1)} \pr_{V}^{*}V\otimes \mc{O}_{\mb{P}_{X}(V)}(-1)\xrightarrow{\taut_{V}} \mc{O}_{\mb{P}_{X}(V)}. 
  \end{equation*}
where $\taut_{V}$ is the tautological morphism on $\mb{P}_{X}(V)$. By Proposition 9.7.9 of \cite{EGA}, we have
\begin{equation}
  \label{eq:projectivization}
\mb{P}_{X}(U)\cong \mc{Z}(\taut_{s})\subset\mb{P}_{X}(V), \text{ and }\mc{O}_{\mb{P}_{X}(U)}(1)=\mc{O}_{\mb{P}_{X}(V)}(1)|_{\mb{P}_{X}(U)}.  
\end{equation} 
We say that $U$ is regular if $\taut_{s}$ is regular on $\mb{P}_{X}(V)$. In this case, the derived pushforward $\mathrm{R}\pr_{U*}(\mathcal{O}_{\mb{P}_{X}(U)}(m))$ was computed in Lemma 3.5 of \cite{zhao2020categorical}: we use the cohomological grading for complexes and associate $C_0$ in the complex $$C:=\{\cdots \to C_{-1}\to C_0\to 0\}$$ with cohomological degree $0$. Given $k\in \mb{Z}_{\geq 0}$, let $S^{k}C$ be the complex $C^{\otimes k}/(v\otimes u-(-1)^{\mathrm{deg}_{v}}(-1)^{\mathrm{deg}_{u}}u\otimes v)$ where $\mathrm{deg}_{u}$ and $\mathrm{deg}_{v}$ are the cohomological degrees of $u$ and $v$, respectively, and $\wedge^{k}C:=(S^{k}(C[1]))[-k]$. Let $r:=\rank(U):=\rank(V)-\rank(W)$. Then
  \begin{equation}
\label{hart}
     \mathrm{R}pr_{U*}(\mathcal{O}_{\mb{P}_{X}(U)}(m))\cong  \begin{cases}
      S^{m}U & m\geq 0\\
      0 & -r<m<0 \\
      \mathrm{det}(U)^{-1}\otimes \wedge^{-m-r}(U^{\vee})[1-r]& m\leq -r  
    \end{cases}
  \end{equation}

\subsection{The conormal sheaf of a closed embedding and intrinsic blow-ups}
Given a scheme $X$, let $\Omega_X$ be the cotangent sheaf of $X$, whose fiber at each closed point $x\in X$ is the dual of the tangent space $T_{X,x}$. Let $f:Z\to Y$ be a closed embedding of schemes with the ideal sheaf $\mc{I}$. The conormal sheaf of $Z$ in $Y$ is defined as $C_{Z/Y}:=\mc{I}/\mc{I}^{2}$; it fits into the right exact sequence 
\begin{equation}
  \label{eq:conormal}
  C_{Z/Y}\to \Omega_{Y}|_{Z}\to \Omega_{Z}\to 0.
\end{equation}
If $Z$ is smooth, then \cref{eq:conormal} is also left exact, since the higher terms of the cotangent complex of $Z$ vanish. In this case, for any closed point $z\in Z$, the normal space $N_{z}:=T_{Y,z}/T_{Z,z}$ is the dual space of $C_{Z/Y}|_{z}$. For closed embeddings of schemes $Z\subset Y\subset X$, we have the right exact sequence of conormal sheaves
\begin{equation*}
  C_{Y/X}|_{Z}\to C_{Z/X}\to C_{Z/Y}\to 0.
\end{equation*}

Let $f:B\to A$ be a morphism of commutative algebras, and let $I$ be an ideal of $B$. Then $I\otimes_{B}A\to IA$ is surjective, and thus the canonical morphism 
\begin{equation*}
  I/I^{2}\otimes_BA\to IA/(IA)^{2}
\end{equation*}
is also surjective. More generally, let $\rho:Y\to Z$ be a closed embedding of schemes. Given a Cartesian square of schemes
  \begin{equation}
    \label{eq:222}
  \begin{tikzcd}
    \mc{R} \arrow[r, hook] \arrow[d,"\sigma"] & \mc{S} \arrow[d, "\theta"]\\
    Z \arrow[r, hook,"\rho"] & Y,
  \end{tikzcd}
\end{equation}
there is a canonical surjective morphism $\sigma^{*}C_{Z/Y}\to C_{\mc{R}/\mc{S}}$. If both $\mc{R}$ and $Z$ are smooth, then for any closed point $r\in \mc{R}$, the Cartesian square \cref{eq:222} induces a map of normal vector spaces 
\begin{equation*}
  N_{r}\to N_{\sigma(r)},
\end{equation*}
which is dual to the morphism $\sigma^{*}C_{Z/Y}\to C_{\mc{R}/\mc{S}}$ at $r$. 

Let $\alpha:Z\subset Y$ be a regular closed embedding of smooth varieties of codimension $2$ with a commutative diagram of exact sequences of coherent sheaves on $Y$:
 \begin{equation}
  \label{eq:aaaa}
  \begin{tikzcd}
    0 \ar{r} & W_0 \ar{r}{w} \ar{d}{u_0} & W_1 \ar{r}{m} \ar{d}{u_1} & B\otimes I_{Z} \ar{d}\ar{r} & 0\\
    0 \ar{r} & V_0 \ar{r}{v} & V_1 \ar{r}{n} & B \ar{r} & 0,
  \end{tikzcd}
\end{equation}
such that $W_0,W_1,V_0,V_1$ are locally free, $B$ is a line bundle and $I_Z$ is the ideal sheaf of $Z$.
Restricting \cref{eq:aaaa} to $Z$, we obtain a commutative diagram of exact sequences
\begin{equation}
  \label{eq:bbbb}
  \begin{tikzcd}
  0  \ar{r} & B\otimes \wedge^2 C_{Z/Y}^{\vee} \ar{d} \ar{r}& \ar{d} W_0|_{Z} \ar{r}{w|_{Z}} \ar{d}{u_0|_{Z}} & W_1|_{Z} \ar{r}{m|_Z} \ar{d}{u_1|_{Z}} & B\otimes C_{Z/Y} \ar{r} \ar{d}{0} & 0\\
    & 0 \ar{r} & V_0|_{Z} \ar{r}{v|_{Z}} & V_1|_{Z} \ar{r}{n|_Z} & B|_{Z} \ar{r} & 0,
  \end{tikzcd}
\end{equation}
which induces a unique morphism $\psi: W_1|_{Z}\to V_0|_{Z}$ that commutes with \cref{eq:bbbb}. The second row of \cref{eq:aaaa} induces a closed embedding $\beta:Y\hookrightarrow \mb{P}_Y(V_1)$ such that $C_{Y/\mb{P}_Y(V_1)}\cong V_0\otimes B^{-1}$. Restricting to $Z$, we obtain a Cartesian diagram
\begin{equation}
  \label{eq:cccc}
  \begin{tikzcd}
    Z \ar{r} \ar{d}{\beta|_Z} & Y \ar{d}{\beta}\\
    \mb{P}_Z(V_1|_Z) \ar{r} & \mb{P}_Y(V_1)
  \end{tikzcd}
\end{equation}
and hence $C_{Z/\mb{P}_Y(V_1)}\cong C_{Z/Y}\oplus (V_0|_Z\otimes B^{-1})$.
\begin{lemma}
\label{lemmaA}
  
  Let $D:=\mb{P}_Y(u_1)\subset \mb{P}_Y(V_1)$ be the projectivization of $u_1:W_1\to V_1$. Then the closed embedding $Z\subset \mb{P}_Y(V_1)$ in \cref{eq:cccc} factors through $D$, and the following sequence is right exact as a sequence of coherent sheaves on $Z$:
  \begin{equation}
    \label{eq:dddd}
    W_1|_{Z}\xrightarrow{\psi|_Z\oplus m|_Z} V_1|_{Z}\oplus B\otimes C_{Z/Y}\xrightarrow{B\otimes \theta} B\otimes C_{Z/D}\to 0,
  \end{equation}
  where $\theta:C_{Z/\mb{P}_Y(V_1)}\to C_{Z/D}$ is the canonical morphism induced by the closed embedding $D\subset \mb{P}_Y(V_1)$. 
\end{lemma}
\begin{proof}
The argument is local; hence we may assume that the second column of \cref{eq:aaaa} splits and 
\begin{equation*}
  u_1\cong u_+\oplus m_+,
\end{equation*}
where $u_+:W_1\to V_0$ and $m_+:W_1\to B$ commute with the diagrams in \cref{eq:aaaa} and $m_+|_{Z}=0$. 
Moreover, after restricting to a smaller open neighborhood, we may assume that $V\cong \mc{O}_Y^k$, $W\cong \mc{O}_Y^l$, and $B\cong \mc{O}_Y$. Then the image of $m_+$ is exactly $I_Z$ and
$\mb{P}_Y(V_1)\cong Y\times \mb{P}^{m}$, where the embedding $\beta$ has homogeneous coordinates $[0:\cdots:0:1]$ in $\mb{P}^m$. Moreover, the tautological morphism $\taut_{V_1}:V_1\to \mc{O}_{\mb{P}_Y(V_1)}(1)$ is $0$ on $Z$, and thus $Z\subset D$. 

To compute the conormal sheaf $C_{Z/D}$, we can restrict to $Y\times \mb{A}^n\subset Y\times \mb{P}^n$ where the last projective coordinate is nonzero. Over the open subset $Y\times \mb{A}^n$, we have
$\mc{O}_{\mb{P}_Y(V_1)}(1)\cong \mc{O}_{Y\times \mb{A}^n}$ and the tautological morphism 
\begin{equation*}
  \taut_{V_1}:\mc{O}_{Y\times \mb{A}^n}^{k+1}\to \mc{O}_{Y\times \mb{A}^n}
\end{equation*}
is given by the vector $(x_1,x_2,\cdots,x_k,1)$, where $x_1,\cdots,x_k$ are the coordinates of $\mb{A}^k$. 
If we represent $u_+$ by an $l\times k$ matrix and $m_+$ by a column vector
\begin{equation*}
  u_+=\begin{pmatrix}
    a_{11} & a_{21} & \cdots & a_{l1} \\
    a_{12} & a_{22} & \cdots & a_{l2} \\
    \vdots & \vdots & \ddots & \vdots \\
    a_{1k} & a_{2k} & \cdots & a_{lk},
  \end{pmatrix}\quad
  m_+=\begin{pmatrix}
    y_1 \\
    y_2 \\
    \vdots \\
    y_l.
  \end{pmatrix}
\end{equation*}
The composition morphism $\taut_{V_1}\circ u_+:\mc{O}^l\to \mc{O}_{Y\times \mb{A}^n}$ is then given by the vector
\begin{equation*}
  (x_1,\cdots,x_k,1)(u_+,m_+).
\end{equation*}
The image of $\taut_{V_1}\circ u_+$, denoted by $J$, is the ideal sheaf of $D$. The ideal sheaf of $Z$ in $Y\times \mb{A}^{n}$ is $H:=I_Z+(x_1,\cdots,x_k)$. Since $J\subset H$, let $\widetilde{taut}:\mc{O}^l\to H$ be the map induced by $\taut_{V_1}\circ u_+$. Then we have 
\begin{equation*}
  \psi|_Z\oplus m|_Z=\widetilde{taut}|_Z:W_1|_Z\to H/H^2.
\end{equation*}
Its image is $(J+H^2)/H^2$. Hence its cokernel is $H/(J+H^2)$, which is exactly the conormal sheaf $C_{Z/D}$.
\end{proof}

Dualizing \cref{eq:bbbb}, we obtain a long exact sequence of coherent sheaves on $Z$:
\begin{equation*}
  0 \to B^{\vee}\otimes C_{Z/Y}^{\vee} \to W_1^{\vee}|_{Z} \to W_0^{\vee}|_{Z} \to B^{\vee}\otimes \wedge^2 C_{Z/Y}^{\vee} \to 0.
\end{equation*}
It induces a closed embedding of $Z$ into $\mb{P}_Z(W_0^{\vee})$ and hence a closed embedding $\lambda:Z\hookrightarrow\mb{P}_Y(W_0^{\vee})$. 

\begin{lemma}
  \label{lemmaB}
  The closed embedding $\lambda:Z\hookrightarrow \mb{P}_Y(W_0^\vee)$ factors through $\mb{P}_Y(u^{\vee}_0)$. Moreover, the conormal sheaf $C_{Z/\mb{P}_Y(u_0^{\vee})}$ is the cokernel of the morphism 
  $$\wedge^2C_{Z/Y}\otimes B\otimes \psi^{\vee}: \wedge^2C_{Z/Y}\otimes B\otimes V_0^{\vee}|_Z\to \wedge^2C_{Z/Y}\otimes B\otimes W_1^{\vee}|_Z.$$
\end{lemma}
\begin{proof}
First, by the definition of $\lambda$, we have $\mc{O}_{\mb{P}_Y(W_0^\vee)}(1)|_Z\cong \wedge^2 C_{Z/Y}\otimes B^{\vee}$. Hence there is a canonical map from $\wedge^2C_{Z/Y}\otimes B\otimes W_1^{\vee}|_Z$ to $C_{Z/\mb{P}_Y(W_0^\vee)}$. We may now reduce to the local case and assume that there is a global morphism $\widetilde{\psi}:W_1\to V_0$ whose restriction to $Z$ is $\psi$. 
To compute the conormal sheaf $C_{Z/\mb{P}_Y(u_0^{\vee})}$, note that over $\mb{P}_Y(W_0^{\vee})$, $\mb{P}_Y(u_0^{\vee})$ and 
$Z$ are the zero locus of the morphisms
\begin{align*}
   \taut_{u_0^{\vee}}:\pr_{W_0}^{*}V_0^{\vee}\otimes \mc{O}_{\mb{P}_Y(W_0^{\vee})}(-1)\to \mc{O}_{\mb{P}_Y(W_0^{\vee})}  \\ \taut_{W_0^{\vee}}\circ \pr_{W_0^{\vee}}^*w:\pr_{W_0}^{*}V_0^{\vee}\otimes \mc{O}_{\mb{P}_Y(W_0^{\vee})}(-1)\to \mc{O}_{\mb{P}_Y(W_0^{\vee})}
\end{align*}
respectively. Hence the conormal sheaf $C_{Z/\mb{P}_Y(W_0^{\vee})}$ is the cokernel of the restriction of the morphism $\pr_{W_0^{\vee}}^*\widetilde{\psi}\otimes \mc{O}_{\mb{P}_Y(W_0^{\vee})}(-1)$ to $Z$, which is exactly the morphism $\wedge^2C_{Z/Y}\otimes B\otimes \psi^{\vee}$.
\end{proof}

\subsection{The intrinsic blow-up}
\label{s27}
Let $\rho:Z\subset Y$ be a closed embedding of schemes, where $Z$ is a smooth variety. We define a new scheme over $Y$, called the \textbf{intrinsic blow-up} of $Y$ along $Z$ and denoted by $\mc{B}l_{Z}Y$, with a projection morphism $pr_{\rho}:\mc{B}l_{Z}Y\to Y$ such that 
\begin{enumerate}[leftmargin=*]
  \item The restriction $pr_{\rho}|_{Y-Z}$ is an isomorphism;
  \item The fiber of $pr_{\rho}$ over $Z$ is the projectivization $\mb{P}_{Z}(C_{Z/Y})$ of the conormal sheaf of $Z$ in $Y$, and $pr_{\rho}|_Z$ is the projection morphism $\mb{P}_{Z}(C_{Z/Y})\to Z$.
  \item For any Cartesian square \cref{eq:222} such that $\mc{R}$ is a Cartier divisor of $\mc{S}$, there exists a unique morphism $\tilde{\theta}:\mc{S}\to \mc{B}l_{Z}Y$ such that $\theta=pr_{\rho}\circ\tilde{\theta}$ and the restriction morphism $\tilde{\theta}|_{\mc{R}}:\mc{R}\to \mb{P}_{Z}(C_{Z/Y})$ is the morphism induced by the surjection $\sigma^{*}C_{Z/Y}\to C_{\mc{R}/\mc{S}}$.
\end{enumerate}

If both $Y$ and $Z$ are smooth, the classical blow-up $\Bl_{Z}Y$ satisfies the above properties. If $Y$ is the zero locus of a cosection $u:V\to \mc{O}_X$ on a smooth variety $X$, let $g:Z\subset X$ be the closed embedding and $pr_g:\Bl_ZX\to X$ be the projection morphism. Let $R:=\mathbb{P}_Z(C_{Z/X})$ be the exceptional divisor of $\Bl_ZX$. The image of 
$$pr_g^*u\otimes \mc{O}_{\Bl_{Z}X}(R):pr_g^*V\otimes \mc{O}_{\Bl_Z{X}}(R)\to \mathcal{O}_{\Bl_{Z}X}(R)$$
 lies in $\mathcal{O}_{\Bl_{Z}X}$ and thus induces a morphism $\tilde{u}:pr_g^*V\otimes \mathcal{O}(R)\to \mathcal{O}$.
 \begin{lemma}
  \label{l27}
  The zero locus of the cosection $\tilde{u}$ satisfies the above three properties. 
 \end{lemma} 
\begin{proof}
  Property (1) is obvious because $\tilde{u}=u$ on $\Bl_{Z}X-R$. For property (3), note that, given a Cartesian diagram \cref{eq:222}, we also have $\mc{R}\cong Z\times_X \mc{S}$. Hence there exists a map $\tilde{g}:\mc{S}\to \Bl_ZX$ such that the map from $\mc{S}$ to $X$ factors through $\tilde{g}$. It remains to prove that $\tilde{g}^*(\tilde{u})=0$, which follows from the fact that on $\mc{S}-\mc{R}$
  \begin{equation*}
    \tilde{g}^*(\tilde{u})=g^*(u)=0.
  \end{equation*}

For property (2), we reduce to the case where $X\cong \mathrm{Spec}(A)$, $V\cong \mc{O}^n$, and $u$ is given by elements $a_1,\cdots,a_n\in \mc{A}$. Let $I$ be the ideal of $Z$, and let $J$ be the ideal generated by $a_1,\cdots, a_n$.  Then $\Bl_ZX=\mathrm{Proj}_A(\bigoplus_{k=0}^{\infty}I^k)$, and the zero locus of $\tilde{u}$ is given by the projective spectrum
\begin{equation*}
 \mathrm{Proj}_A(\bigoplus_{k=0}^{\infty}I^k/I^{k-1}J)=\mathrm{Proj}_A(\bigoplus_{k=0}^{\infty}\mathrm{coker}(I^{k-1}\otimes V\to I^k)).
\end{equation*}
Hence its fiber over $Z$ is given by
\begin{equation*}
  \mathrm{Proj}_{A/I}(\mathrm{coker}(\bigoplus_{k=0}^{\infty}(I^{k-1}/I^{k-2})^{\oplus n}\to I^{k}/I^{k-1}))
\end{equation*}
which is the projectivization of $C_{Z/Y}$. 
\end{proof}

\begin{remark}
  The construction of $\mc{B}l_ZY$ above relies on an embedding of $Y$ into an ambient smooth variety $X$. We do not prove the uniqueness of this construction in this paper, leaving it to future work. 
\end{remark}
In this paper, we denote the zero locus of the cosection $\tilde{u}$ by $\mc{B}l_ZY$ and let $pr_\rho:\mc{B}l_ZY\to Y$ be the projection morphism. In the setting of \cref{l27}, we call the restriction of $\mc{O}(R)$ to $\mc{B}l_ZY$ the exceptional divisor of $\mc{B}l_ZY$. 
\begin{lemma}
  \label{l28}
  Let $d:=\dim(X)-\rank(V)$. 
  If $Y-Z$ is a $d$-dimensional smooth variety and $\mb{P}_Z(C_{Z/Y})$ is a $(d-1)$-dimensional smooth variety, then $\mc{B}l_ZY$ is a smooth variety of dimension $d$. Moreover, if $\dim(Z)<d$, then $Y$ is a locally complete intersection variety, and we have 
  \begin{equation}
    \label{eq:pushforward}
    \mathrm{R}pr_{\rho*}(\mathcal{O}_{\mb{B}l_ZY})\cong \mathcal{O}_Y.
  \end{equation}
\end{lemma}
\begin{proof}
  We note that $\mc{B}l_ZY$ is the zero locus of the cosection $\tilde{u}$ on $\Bl_ZX$, and $\mb{P}_Z(C_{Z/Y})$ is the zero locus of $\tilde{u}$ on $R$. Hence, to prove the smoothness of $\mc{B}l_ZY$, it suffices to use the following local argument: let $A$ be a local ring with maximal ideal $m$, and let $f\in m$ be a non-zero divisor. If $A/f$ is a regular local ring with $\dim(A/f)=\dim(A)-1$, then $A$ is also a regular local ring. The proof is similar to that of Krull's principal ideal theorem, so we omit the details. If $\dim(Z)< \dim(Y-Z)$, then $Y$ is a locally complete intersection variety by dimension counting, and hence the Koszul complex $\wedge^{\bullet}f$ is a resolution of $\mathcal{O}_Y$. 

  Finally, we prove \cref{eq:pushforward}. First, consider the case where $Y$ is smooth. Then \cref{eq:pushforward} is well known; moreover, by the resolutions of exceptional divisors, we have 
    \begin{equation*}
    \mathrm{R}pr_{\rho*}(\mathcal{O}_{\mb{B}l_ZY}(kR))\cong \mathcal{O}_Y.
  \end{equation*}
  when $\dim(Z)-\dim(Y)< k\leq 0$. For the general case, note that $\mc{O}_{\mc{B}l_ZY}\cong \wedge^{\bullet}\tilde{u}$ and hence \cref{eq:pushforward} follows by comparing the pushforward of $\wedge^{\bullet}\tilde{u}$ in each degree of the complex.
\end{proof}

\subsection{Fourier--Mukai transforms associated to an algebraic surface}
Given a scheme $X$, let $\D^{b}_{\mathrm{coh}}(X)$ be the bounded derived category of coherent sheaves on $X$.

\begin{definition}[Definition 6.1 of \cite{zhao2020categorical}]
  \label{FMS}
   Let $S_1$ and $S_2$ be two copies of $S$, used to distinguish the two factors of $S\times S$. Given smooth varieties $X,Y,Z$ and Fourier--Mukai kernels $P\in \D^{b}_{\mathrm{coh}}(X\times Y\times S_1)$ and $Q\in \D^{b}_{\mathrm{coh}}(Y\times Z\times S_{2})$, we define the composition $QP\in \D^{b}_{\mathrm{coh}}(X\times Z\times S_1\times S_2)$ by
$$QP:=\mathrm{R}\pi_{13_{*}}(L\pi_{12}^*P\otimes L\pi_{23}^*Q),$$
where $\pi_{12},\pi_{23}$, and $\pi_{13}$ are the projections from $X\times Y\times Z\times S_1\times S_2$ to $X\times Y\times S_1$, $Y\times Z\times S_2$, and $X\times Z\times S_1\times S_2$, respectively.
\end{definition}

\subsection{Derived algebraic geometry}
\label{sec5}
We recall some facts about dg-schemes and derived schemes that will only be used in \cref{sec6}.
A dg-scheme consists of a scheme $X_{0}$ and quasi-coherent sheaves $\mc{O}_{X}:=\{\mc{O}_{X,i}\}_{i\geq 0}$ on $X_{0}$ such that $\mc{O}_{X,0}=\mc{O}_{X_{0}}$, together with a cdga structure, i.e., maps $\delta:\mc{O}_{X,i}\to \mc{O}_{X,i-1}$ and $\bullet:\mc{O}_{X,i}\otimes \mc{O}_{X,j}\to \mc{O}_{X,i+j}$ compatible with the differentials and multiplication. We denote by $\pi_{0}X$ the closed subscheme of $X_{0}$ such that $\mc{O}_{\pi_{0}X}=H^{0}(\mc{O}_{X})$.

A derived scheme $X$ consists of a pair $(X,\mc{O}_{X})$, where $X$ is a topological space and $\mc{O}_{X}$ is a sheaf of commutative simplicial rings on $X$ such that the ringed space $(X,\pi_{0}(\mc{O}_{X}))$ is a scheme and the homotopy sheaf $\pi_{i}(\mc{O}_{X})$ is a quasi-coherent sheaf over $(X,\pi_{0}(\mc{O}_{X}))$ for all $i>0$. Given a dg-scheme $(X_{0}, \mc{O}_{X})$, there is a canonical closed embedding $i:\pi_{0}X\to X_{0}$, and $(\pi_{0}X,i^{-1}\mc{O}_{X})$ is a derived scheme.

Given any scheme $X$, $(X,\mc{O}_{X})$ is also a derived scheme, which we call a classical scheme in this paper. Let $\mathrm{dSch}$ be the $\infty$-category of derived schemes. For any two morphisms 
  \begin{equation}
    \label{eq:cartesian}
    f:X\to Z, g:Y\to Z
  \end{equation}
in $\mathrm{dSch}$, \cite{toen2008homotopical} shows that the Cartesian product in $\mathrm{dSch}$ exists; we denote it by $X\times_{Z}^{\mb{L}}Y$.
\begin{lemma}
  \label{n610}
  If $f$ and $g$ in \cref{eq:cartesian} are regular embeddings of smooth classical schemes, and $\dim(X\times_{Z}Y)=\dim(X)+\dim(Y)-\dim(Z)$, then the Cartesian product $X\times_{Z}^{\mb{L}}Y$ is isomorphic to the classical fiber product $X\times_{Z}Y$.
\end{lemma}
\begin{proof}
  It follows from Proposition 3.2 of \cite{neguct2018hecke} that $\mathrm{tor}^{i}_{Z}(\mc{O}_{X},\mc{O}_{Y})=0$ when $i>0$.
\end{proof}

\begin{example}[Derived zero locus of a cosection]
  Let $X$ be a scheme, and let $f:V\to \mc{O}_X$ be a morphism of locally free sheaves. The Koszul complex $\wedge^{\bullet}f$ defines a cdga structure on $X$ and thus induces a dg-scheme, which we denote by $\mb{R}\mc{Z}(f)$. Let $\mathrm{Taut}(V^{\vee})$ be the tautological bundle corresponding to $V^{\vee}$. Regarded as a derived scheme, $\mb{R}\mc{Z}(f)$ is the derived Cartesian product of $f^{\vee}:X\to \mathrm{Taut}(V^{\vee})$ and the zero section $0:X\to \mathrm{Taut}(V^{\vee})$.
\end{example}

\begin{definition}[The bounded derived category]
  Given a derived scheme $\mf{X}$, we define the dg-category
\begin{equation}
  \label{lim}
  \mathrm{L}_{qcoh}(\mf{X}):=\lim_{\mf{U}\to \mf{X}}\mathrm{L}_{qcoh}(\mf{U}).
\end{equation}
Here $\mf{U}=\mathrm{Spec} A$ is an affine derived scheme for a cdga $A$, and the category $\mathrm{L}_{qcoh}(\mf{U})$ is the dg-category of dg-modules over $A$ localized by quasi-isomorphisms. The homotopy category of $\mathrm{L}_{qcoh}(\mf{X})$ contains a triangulated subcategory
$\D^{b}_{\mathrm{coh}}(\mf{X})$ consisting of objects with bounded coherent cohomology.

\end{definition}
\begin{example}
  \label{m57}
  Given a dg-scheme $(X_{0},\mc{O}_{X})$, let $\mf{X}$ be the derived scheme $(\pi_{0}X, i^{-1}\mc{O}_{X})$, where $i:\pi_{0}X\to X$ is the canonical embedding. Given a dg-$\mc{O}_{X}$-module $\mc{F}$ whose cohomology sheaves are bounded and coherent $\mc{O}_{\pi_{0}X}$-modules, $i^{-1}\mc{F}\in D^{b}_{\mathrm{coh}}(\mf{X})$. Given dg-$\mc{O}_{X}$-modules $\mc{F},\mc{G},\mc{H}$ whose cohomology sheaves are bounded and coherent $\mc{O}_{\pi_{0}X}$-modules, together with $\mc{O}_{X}$-morphisms
  \begin{equation}
    \label{exact}
 0\to \mc{F}\xrightarrow{a}\mc{G}\xrightarrow{b} \mc{H}\to0,
  \end{equation}
  that are exact in each degree, \cref{exact} can be lifted to a triangle in $\D^{b}_{\mathrm{coh}}(\mf{X})$.
\end{example}

\begin{definition}[Quasi-smooth morphisms]
  A morphism of derived schemes $f:\mf{X}\to \mf{W}$ is called quasi-smooth if the $f$-relative cotangent complex $\mb{L}_{f}$ is perfect and, for any point $x\in \pi^{0}\mf{X}$, the restriction $\mb{L}_{f}|_{x}$ is of cohomological amplitude $[-1,1]$.
\end{definition}

Morphisms between classical smooth varieties are always l.c.i. and thus quasi-smooth. Finally, we introduce the pullback and pushforward functors:
\begin{lemma}[Section 4.2 of \cite{porta2019categorification} and Proposition 1.4 of \cite{toen2012proper}]
  \label{n613}
Given a morphism of derived schemes $f:X\to Y$, there is a natural pullback functor $Lf^{*}:\D_{coh}^{b}(Y)\to \D_{coh}^{b}(X)$ if $f$ is quasi-smooth and a natural pushforward functor $Rf_{*}:\D_{coh}^{b}(X)\to \D_{coh}^{b}(Y)$ if $f$ is proper. Moreover, for any derived square
    \begin{equation*}
      \begin{tikzcd}
        X' \ar{r}{g}\ar{d}{q} & X \ar{d}{p}  \\
        Y' \ar{r}{f} & Y
      \end{tikzcd}
    \end{equation*}
  such that all derived schemes are quasi-compact and quasi-separated, $f$ is proper, and $p$ is quasi-smooth, there is a canonical isomorphism of functors:
  $$f^{*}p_{*}\Rightarrow q_{*}g^{*}: \D^{b}_{coh}(X)\to \D^{b}_{coh}(Y').$$
\end{lemma}

\section{Nested, triple, and quadruple moduli spaces of sheaves}
\label{sec:n}
Let $S$ be a smooth projective algebraic surface over $\mb{C}$ with an ample line bundle $H$. Let $q:=K_{S}$ be the canonical line bundle of $S$. Fix $(r,c_{1})\in \mb{N}^{>0}\times H^{2}(S,\mb{Z})$, and let $\mc{M}_k$ be the moduli space of Gieseker $H$-stable sheaves $\mc{F}$ such that $\rank(\mc{F})=r$, $c_{1}(\mc{F})=c_{1}$, and $c_{2}(\mc{F})=k$ for any integer $k$. We assume
 \begin{equation*}
  \gcd(r, c_{1}\cdot H) = 1,  \qquad  q \cong \mc{O}_S \text{ or } c_1(q) \cdot H < 0.
\end{equation*}

Now we recall the geometry of the nested, triple, and quadruple moduli spaces of sheaves from \cite{negut2017shuffle,neguct2018hecke} and give a new description of the quadruple moduli spaces via intrinsic blow-ups in \cref{p36}. Throughout this paper, $\mathrm{const}$ denotes the constant in \cref{sec:1.2}.
\subsection{Nested moduli space $\mf{Z}_{k,k+1}$} Given an integer $k$, we recall the definition of the nested moduli space of sheaves $\mf{Z}_{k,k+1}$ in \cref{e11}:
\begin{equation*}
   \mf{Z}_{k,k+1}:=\{(\mc{F}_{-1}\subset_{x} \mc{F}_{0})|\mc{F}_{-1}, \mc{F}_{0} \text{ are stable }, c_{2}(\mc{F}_{0})=k\},
\end{equation*}
There is a tautological line bundle $\mc{L}$ on $\mf{Z}_{k,k+1}$ whose fiber at a closed point is $\mc{F}_0/\mc{F}_{-1}$. The projections to $\mc{F}_{0},\mc{F}_{-1}$, and $x$, respectively, induce projection morphisms
$$p_{k}:\mf{Z}_{k,k+1}\to \mc{M}_{k}, \quad q_{k}:\mf{Z}_{k,k+1}\to \mc{M}_{k+1},\quad \pi_{k}:\mf{Z}_{k,k+1}\to S.$$
Let $\mc{U}_k$ be the universal sheaf on $\mc{M}_{k}\times S$. By Proposition 2.14 of \cite{neguct2018hecke}, there exists a resolution of the universal sheaf $\mc{U}_k$ by locally free sheaves:
\begin{equation*}
        0 \to \mc{W}_{k}\xrightarrow{\psi_{k}}\mc{V}_{k}\to \mc{U}_{k}\to 0.
  \end{equation*}
By Propositions 2.19 and 2.26 of \cite{neguct2018hecke}, if we abuse notation and write $\mc{U}_{k}:=\{\mc{W}_{k}\xrightarrow{\psi_{k}}\mc{V}_{k}\}$, then 
\begin{equation}
  \label{p25}
    \mf{Z}_{k,k+1}\cong \mb{P}_{\mc{M}_{k}\times S}(\mc{U}_{k}),\quad \mf{Z}_{k-1,k}\cong \mb{P}_{\mc{M}_{k}\times S}(q\mc{U}_{k}^{\vee}[-1]).
\end{equation}
where $p_{k}\times \pi_{k}$ and $q_{k}\times \pi_{k}$ are the respective projection morphisms. Moreover, the projectivizations in \cref{p25} are regular, and we have 
 $$\mc{L}\cong \mc{O}_{\mb{P}_{\mc{M}_{k}\times S}(\mc{U}_{k})}(1)\cong \mc{O}_{\mb{P}_{\mc{M}_{k}\times S}(q\mc{U}_{k}^{\vee}[-1])}(-1).$$
Equation \cref{hart} implies the following formula: 
  \begin{equation}
      \label{cor:5.3ext}
     R(p_{k}\times \pi_{k})_{*}(\mc{L}^{m})=
     \begin{cases}
       S^{m}(\mc{U}_{k}) & m\geq 0 \\
       0 & -r<m<0 \\
     \det(\mc{U}_{k})^{-1}\wedge^{-m-r}(\mc{U}_{k}^{\vee})[1-r] & m\leq -r
     \end{cases}
   \end{equation}

Over $\mf{Z}_{k,k+1}\times S$, we abuse notation and let $\mc{U}_{k}$ and $\mc{U}_{k+1}$ denote the universal sheaves, respectively. Then $\mf{Z}_{k,k+1}$ is a closed subscheme of $\mf{Z}_{k,k+1}\times S$ via the graph of $\pi_k$. By Section 2.20 of \cite{neguct2018hecke}, over $\mf{Z}_{k,k+1}$ there exists a commutative diagram of exact sequences of sheaves:
  \begin{equation}
  \label{eq:1.7}
  \begin{tikzcd}
    0\ar{r} & \mathcal{W}_{k+1}\ar{r}{\mf{w}_{k,k+1}}\ar{d}{\psi_{k+1}} & \mathcal{W}_{k} \ar{r}\ar{d}{\psi_{k}} &  \mc{B}_{k,k+1}\otimes \mathcal{I}_{\mf{Z}_{k,k+1}} \ar{r}\ar{d} & 0 \\
    0\ar{r} & \mathcal{V}_{k+1}\ar{r}{\mf{v}_{k,k+1}} & \mathcal{V}_{k} \ar{r} & \mc{B}_{k,k+1} \ar{r} & 0
  \end{tikzcd}
\end{equation}
which matches the toy model \cref{eq:aaaa}, i.e.,
\begin{enumerate}
  \item $\mc{W}_{k},\mc{W}_{k+1},\mc{V}_{k},\mc{V}_{k+1}$ are locally free sheaves on $\mf{Z}_{k,k+1}\times S$;
  \item $\mc{U}_{k}=\mathrm{coker}(\psi_{k})$ and $\mc{U}_{k+1}=\mathrm{coker}(\psi_{k+1})$ and $\psi_{k}$ and $\psi_{k+1}$ are injective morphisms;
  \item $\mc{B}_{k,k+1}$ is a line bundle on $\mf{Z}_{k,k+1}\times S$ such that the restriction $\mc{B}_{k,k+1}|_{\mf{Z}_{k,k+1}}$ is isomorphic to $\mc{L}$, and $\mc{I}_{\mf{Z}_{k,k+1}}$ is the ideal sheaf of $\mf{Z}_{k,k+1}$ in $\mf{Z}_{k,k+1}\times S$.
\end{enumerate}
Restricting \cref{eq:1.7} to $\mf{Z}_{k,k+1}$, we obtain the following commutative diagram
 \begin{equation}
   \label{commute}
  \begin{tikzcd}
  0\ar{r}&  q\mc{L}\ar{d}\ar{r} & \mc{W}_{k+1} \ar{r}{\mf{w}_{k,k+1}} \ar{d}{\psi_{k+1}}& \mc{W}_{k}\ar{r}{\mf{w}_{k}} \ar{d}{\psi_{k}} & \mc{L}\otimes \Omega_S \ar{r} \ar{d}{0} & 0  \\
    & 0 \ar{r} & \mc{V}_{k+1} \ar{r}{\mf{v}_{k,k+1}} & \mc{V}_{k}\ar{r} & \mc{L} \ar{r}& 0.
  \end{tikzcd}
\end{equation}
where all rows are long exact sequences. Hence there is a unique morphism 
\begin{equation}
  \label{psik}
  \psi_{k,k+1}:\mc{W}_{k}\to \mc{V}_{k+1}
\end{equation}
which makes \cref{commute} commute.

\subsection{Nested moduli space $\mf{Z}_{k-1,k,k+1}^{\bullet}$} We define $\mf{Z}_{k-1,k,k+1}^{\bullet}$ as the moduli space of triples of sheaves with the same support point:
\begin{align*}
 \mf{Z}_{k-1,k,k+1}^{\bullet}=\{(\mc{F}_{-1}\subset_{x} \mc{F}_{0}\subset_{x} \mc{F}_{1})| \mc{F}_{-1},\mc{F}_{0},\mc{F}_{1}\text{ are stable}, c_{2}(\mc{F}_{0})=k\}.
\end{align*}
Forgetting $\mc{F}_{-1}$ and $\mc{F}_{1}$, respectively, induces morphisms
$$p_{k-1,k}:\mf{Z}_{k-1,k,k+1}^{\bullet}\to \mf{Z}_{k-1,k},\quad q_{k,k+1}:\mf{Z}_{k-1,k,k+1}^{\bullet}\to \mf{Z}_{k,k+1}.$$
Let $\mathcal{L}_{1}$ and $\mathcal{L}_{2}$ be tautological line bundles on $\mf{Z}_{k-1,k,k+1}^{\bullet}$ whose fibers at a closed point are $\mathcal{F}_{0}/\mc{F}_{-1}$ and $\mathcal{F}_{1}/\mathcal{F}_{0}$, respectively. Composing with the projection morphisms gives 
\begin{equation}
  \label{hart1}
  (p_{k-1}\times \pi_{k-1})\circ p_{k-1,k}=(p_{k}\times \pi_{k})\circ q_{k,k+1}.
\end{equation}

Recalling $\psi_{k,k+1}$ from \cref{psik}, by Propositions 2.19 and 2.27 of \cite{neguct2018hecke},
 the scheme $\mf{Z}_{k-1,k,k+1}^{\bullet}$ is smooth of dimension $\mathrm{const}+2rk+1$ and
\begin{equation}
   \label{lemma:n39}
    \mf{Z}_{k-1,k,k+1}^{\bullet}\cong \mb{P}_{\mf{Z}_{k-1,k}}(\psi_{k-1,k}\oplus \mf{w}_{k-1})\cong \mb{P}_{\mf{Z}_{k,k+1}}(q\mc{L}\otimes \psi_{k,k+1}^{\vee}),
\end{equation}
where
  $$\psi_{k-1,k}\oplus \mf{w}_{k-1}:\mc{W}_{k-1}\to \mc{V}_{k}\oplus \mc{L}\otimes \Omega_{S}, \quad q\mc{L}\otimes \psi_{k,k+1}^{\vee}:q\mc{L}\otimes \mc{V}_{k+1}^{\vee}\to q\mc{L}\otimes \mc{W}_{k}^{\vee} $$
  are morphisms of locally free sheaves on $\mf{Z}_{k-1,k}$ and $\mf{Z}_{k,k+1}$, respectively, defined in \cref{commute}. Both projectivizations are regular, and $p_{k-1,k}$ and $q_{k,k+1}$ are the respective projection morphisms. Moreover, we have
    $$\mc{O}_{\mb{P}_{\mf{Z}_{k-1,k}}(\psi_{k-1,k}\oplus \mf{w}_{k-1})}(1)\cong \mc{L}_{1},\quad \mc{O}_{\mb{P}_{\mf{Z}_{k,k+1}}(q\mc{L}\psi_{k,k+1}^{\vee})}(1)\cong \mc{L}_{1}\mc{L}_{2}^{-1}.$$
By \cref{hart1}, \cref{lemma:n39}, and \cref{cor:5.3ext}, we have
  \begin{equation}
   \label{n510} R((p_{k}\times \pi_{k})\circ q_{k,k+1})_{*}\mc{L}_{1}^{i}\mc{L}_{2}^{j}=0, \qquad \text{if } -r\leq i<0.
  \end{equation}
\subsection{Triple and quadruple moduli spaces} \label{s13}
We consider the triple moduli spaces of sheaves whose closed points are given by triples of sheaves
  \begin{align*}
   &  \mf{Z}_{-}=\{(\mc{F}_{0}\subset_{y} \mc{F}_{1}, \mc{F}_{0}'\subset_{x} \mc{F}_{1})|\mc{F}_{0},\mc{F}_{0}',\mc{F}_{1} \text{ are stable, } c_{2}(\mc{F}_{0})=k\} \\
 &  \mf{Z}_{+}=\{(\mc{F}_{0}\supset_{x} \mc{F}_{-1}, \mc{F}_{0}'\supset_{y} \mc{F}_{-1})|\mc{F}_{0},\mc{F}_{0}',\mc{F}_{-1} \text{ are stable, } c_{2}(\mc{F}_{0})=k\}
\end{align*}
There are tautological line bundles $\mathcal{L}_{1}$ and $\mathcal{L}_1'$ over $\mf{Z}_+$ whose fibers at a closed point are $\mc{F}_0/\mc{F}_{-1}$ and $\mc{F}_0'/\mc{F}_{-1}$, respectively, and tautological line bundles $\mc{L}_2$ and $\mc{L}_2'$ over $\mf{Z}_-$ whose fibers at a closed point are $\mc{F}_{1}/\mc{F}_0$ and $\mc{F}_1/\mc{F}_0'$, respectively. There are canonical closed embeddings 
$$\Delta_{-}:\mf{Z}_{k-1,k}\to \mf{Z}_{-}, \quad \Delta_{+}:\mf{Z}_{k,k+1}\to \mf{Z}_{+}$$
whose images consist of triples such that $\mc{F}_{0}\cong \mc{F}_{0}'$. By \cref{p25}, we have 
\begin{equation}
  \label{eq:1.6}
  \mf{Z}_{-}\cong \mb{P}_{\mf{Z}_{k-1,k}\times S}(\mc{U}_{k-1}),\quad \mf{Z}_{+}\cong \mb{P}_{\mf{Z}_{k,k+1}\times S}(q\mc{U}_{k+1}^{\vee}[-1]).
\end{equation}
We now recall the Cartesian diagrams \cref{eq:aaaa} and \cref{eq:1.7}. The conormal sheaf of the graph of $\pi_k$ is exactly $q$. By \cref{lemmaA}, \cref{lemmaB}, and \cref{eq:1.7}, we have 
\begin{equation*}
  C_{\Delta_{+}}\cong \mathrm{coker}(q\mc{L}\otimes \psi_{k,k+1}^{\vee}),\quad C_{\Delta_{-}}\cong \mathrm{coker}(\mc{L}^{-1}(\psi_{k-1,k}\oplus \mf{w}_{k-1})).
\end{equation*}
where $\psi_{k-1,k}\oplus \mf{w}_{k-1}:\mc{W}_k\to \mc{V}_{k+1}\oplus \mc{L}\otimes \Omega_S$ is the morphism defined in \cref{commute}. By \cref{lemma:n39}, we have 
\begin{equation}
  \label{eq:1.8}
  \mb{P}_{\mf{Z}_{k,k+1}}(C_{\Delta_{+}})\cong \mb{P}_{\mf{Z}_{k-1,k}}(C_{\Delta_{-}})\cong \mf{Z}_{k-1,k,k+1}^{\bullet}.
\end{equation}

The quadruple moduli space $\mf{Y}$ consists of quadruples of sheaves:
\begin{equation*}
   \mf{Y}=\{(\mc{F}_{-1}\subset_{x}\mc{F}_{0}\subset_{y} \mc{F}_{1}, \mc{F}_{-1}\subset_{y}\mc{F}_{0}'\subset_{x} \mc{F}_{1})|\mc{F}_{-1},\mc{F}_{0},\mc{F}_{0}',\mc{F}_{1} \text{ are stable, } c_{2}(\mc{F}_{0})=k\}.
\end{equation*}
The subspace $\mf{Z}_{k-1,k,k+1}^{\bullet}\subset\mf{Y}$ consists of quadruples such that $\mc{F}_{0}\cong \mc{F}_{0}'$; let $\Delta_{\mf{Y}}$ denote this closed embedding.
By Proposition 2.39 of \cite{neguct2018hecke}, $\mf{Y}$ is a smooth variety and $\mf{Z}_{k-1,k,k+1}^{\bullet}$ is a smooth divisor of $\mf{Y}$. Moreover, by Proposition 2.41 of \cite{neguct2018hecke} and Claim 3.8 of \cite{negut2017shuffle}, the morphisms
\begin{equation*}
  \alpha_{+}:\mf{Y}\to \mf{Z}_{+}, \quad \alpha_{-}:\mf{Y}\to \mf{Z}_{-}
\end{equation*}
induced by forgetting $\mc{F}_{-1}$ and $\mc{F}_{1}$, respectively, are isomorphisms when restricted to $\mf{Y}-\Delta_{\mf{Y}}$. Let $\mathcal{L}_{1},\mathcal{L}_{2},\mathcal{L}_{1}',\mathcal{L}_2'$ be tautological line bundles over $\mf{Y}$ whose fibers at a closed point are $\mc{F}_{0}/\mc{F}_{-1},\mc{F}_{1}/\mc{F}_{0},\mc{F}_{0}'/\mc{F}_{-1},\mc{F}_{1}/\mc{F}_{0}'$, respectively. Forgetting $\mc{F}_{-1}$ and $\mc{F}_{1}$, respectively, induces the morphisms $\alpha_{-}:\mf{Y}\to \mf{Z}_{-}$ and $\alpha_{+}:\mf{Y}\to \mf{Z}_{+}$. Moreover, we have the following relation between tautological line bundles:
    \begin{equation}
      \label{neq35}
    \mc{L}_{1}^{a}\mc{L}_{1}'^{b}\mc{L}_{2}^{c}\mc{L}_{2}'^{d}=\mc{L}_{1}^{a+d}\mc{L}_{1}'^{b+c}\mc{O}((c+d)\Delta_{\mf{Y}}).
    \end{equation}
Given resolutions of $\mc{U}_{k}$ and $\mc{U}_{k+1}$ on $\mf{Z}_{k,k+1}\times S$ and of $\mc{U}_{k}$ and $\mc{U}_{k-1}$ on $\mf{Z}_{k-1,k}\times S$ as in \cref{eq:1.7}, by \cref{eq:1.6} we have $\mf{Z}_{+}\subset \mb{P}_{\mf{Z}_{k,k+1}\times S}(q\mc{W}_{k+1}^{\vee})$ and $\mf{Z}_{-}\subset \mb{P}_{\mf{Z}_{k-1,k}\times S}(\mc{V}_{k-1})$ as the zero loci of the tautological cosections. By \cref{s27}, the morphisms $\alpha_{+}$ and $\alpha_{-}$ can be lifted to morphisms
   $$\overline{\alpha}_{-}:\mf{Y}\to \mc{B}l_{\mf{Z}_{k-1,k}}\mf{Z}_{-}, \qquad\overline{\alpha}_{+}:\mf{Y}\to \mc{B}l_{\mf{Z}_{k,k+1}}\mf{Z}_{+}.$$

The proof of the following theorem is the key point of this paper and is deferred to \cref{sec:proof}.  
\begin{theorem}
  \label{p36}
  The morphisms $\overline{\alpha}_{+}$ and $\overline{\alpha}_{-}$ are isomorphisms.
\end{theorem}

Consider the following morphisms
\begin{align*}
\theta_{-}:\mf{Z}_{-}\to \mc{M}_{k}\times \mc{M}_{k}\times S\times S, \quad 
(\mc{F}_{0}\supset_{x} \mc{F}_{-1},\mc{F}_{0}'\supset_{y} \mc{F}_{-1})&\to (\mc{F}_{0},\mc{F}_{0}',x,y), \\
\theta_{+}:\mf{Y}\to \mc{M}_{k}\times \mc{M}_{k}\times S\times S, \quad
(\mc{F}_{-1}\subset_{x}\mc{F}_{0}\subset_{y}\mc{F}_{1}, \mc{F}_{-1}\subset_{y}\mc{F}_{0}'\subset_{x}\mc{F}_{1})&\to (\mc{F}_{0},\mc{F}_{0}',y,x).
\end{align*}
where $\theta_{+}=\iota\circ \theta_{-}\circ \alpha_{-}$.
Since 
$$\dim(\mf{Y}_{k-1,k})=\mathrm{const}+2rk+1-r<\mathrm{const}+2rk+2=\dim(\mf{Y}-\Delta_{\mf{Y}}),$$
the projectivization $\mf{Z}_{-}=\mb{P}_{\mf{Z}_{k-1,k}\times S}(\mc{U}_{k-1})$ is regular. By \cref{l28} and \cref{p36}, $R\alpha_{-*}\mc{O}_{\mf{Y}}\cong \mc{O}_{\mf{Z}_{-}}$. Thus we obtain a formula for $R\iota_{*}e_{j}f_{i}$:
\begin{align}
   \label{n42}  R\iota_{*}e_{j}f_{i}&\cong R\iota_{*}R\theta_{-*}\mc{L}_{2}^{j}\mc{L}_{2}'^{i-r}\cong R\iota_{*}R\theta_{-*}R\alpha_{-*}\mc{L}_{2}^{j}\mc{L}_{2}'^{i-r} & \\
                      &\cong R\theta_{+*}\mc{L}_{1}^{i-r}\mc{L}_{1}'^{j}\mc{O}((i+j-r)\Delta_{\mf{Y}})  & \text{ by } \nonumber \cref{neq35} 
\end{align}

 \subsection{The tangent space of the nested moduli space}
 \label{sec:proof}
 To prove \cref{p36}, we first recall the description of the tangent space at a closed point in the nested moduli space of sheaves. Let $d\pi$ be the morphism of tangent spaces induced by $\pi_{k}:\mf{Z}_{k,k+1}\to S$. It is induced as follows: given a closed point $(\mc{F}_{-1}\subset_{x} \mc{F}_{0})\in \mf{Z}_{k,k+1}$, by equation (6.23) of \cite{neguct2018hecke}, the tangent space at this point is given by pairs 
  \begin{equation*}
    \{(w_{-1},w_{0})\in \mathrm{Ext}^1(\mc{F}_{-1},\mc{F}_{-1})\oplus \mathrm{Ext}^1(\mc{F}_{0},\mc{F}_{0})\}
  \end{equation*}
  such that $w_{-1}$ and $w_{0}$ map to the same element in $\mathrm{Ext}^1(\mc{F}_{-1},\mc{F}_{0})$ via the maps induced by the inclusion $\mc{F}_{-1}\subset_x \mc{F}_{0}$. This condition can also be explained by a commutative diagram of exact sequences of sheaves: given a pair $(w_{-1},w_0)\in \mathrm{Ext}^1(\mc{F}_{-1},\mc{F}_{-1})\oplus \mathrm{Ext}^1(\mc{F}_{0},\mc{F}_{0})$, we consider two short exact sequences of coherent sheaves
\begin{align*}
  0 \to \mc{F}_{-1} \to \mc{A} \to \mc{F}_{-1} \to 0, \quad
  0 \to \mc{F}_0 \to \mc{B} \to \mc{F}_0 \to 0,
\end{align*}
whose extension classes are $w_{-1}$ and $w_0$, respectively. Then $w_{-1}$ and $w_0$ map to the same element in $\mathrm{Ext}^1(\mc{F}_{-1},\mc{F}_{0})$ if and only if there exists a commutative diagram of exact sequences of sheaves:
\begin{equation}
  \label{eq:1.9}
  \begin{tikzcd}
    0 \ar{r} & \mc{F}_{-1} \ar{r} \ar[d,hook] & \mc{A} \ar{r} \ar[d,hook] & \mc{F}_{-1} \ar{r} \ar[d,hook] & 0 \\
    0 \ar{r} & \mc{F}_{0} \ar{r} & \mc{B} \ar{r} & \mc{F}_{0} \ar{r} & 0.
  \end{tikzcd}
\end{equation}
Moreover, the quotient of the two exact sequences in \cref{eq:1.9} induces a vector in $\mathrm{Ext}^1(\mb{C}_x,\mb{C}_x)\cong \mathrm{T}_x S$, which is the image of $(w_{-1},w_0)$ under the map of tangent spaces induced by $\pi_k$, denoted by $d\pi$. 

By \cite{neguct2018hecke}, the tangent space at a closed point $(\mc{F}_{-1}\subset_x \mc{F}_{0}\subset_x \mc{F}_1)\in \mf{Z}_{k-1,k,k+1}^{\bullet}$
consists of triples 
  \begin{equation*}
    (w_{-1},w_{0},w_{1})\in \mathrm{Ext}^1(\mc{F}_{-1},\mc{F}_{-1})\oplus \mathrm{Ext}^1(\mc{F}_{0},\mc{F}_{0})\oplus \mathrm{Ext}^1(\mc{F}_{1},\mc{F}_{1})
  \end{equation*}
  such that $w_{-1}$ and $w_{0}$ (resp. $w_0$ and $w_1$) map to the same element in $\mathrm{Ext}^1(\mc{F}_{-1},\mc{F}_{0})$ (resp. $\mathrm{Ext}^1(\mc{F}_{0},\mc{F}_{1})$) and $d\pi(w_0,w_1)=d\pi(w_{-1},w_0)$ in $\mathrm{T}_x S$. On the other hand, we consider the quotient maps $\mc{F}_1\to \mb{C}_x$ and $\mc{F}_{0}\to \mb{C}_x$, and let $\mc{V}_0$ be the cokernel of the map $\mc{F}_0\oplus \mc{F}_1\to \mb{C}_x$. These data induce the following commutative diagram of exact sequences of sheaves:
\begin{equation}
  \label{eq:1.10}
  \begin{tikzcd}
    & 0 \ar{d} & 0 \ar{d} & 0 \ar{d} \\ 
    0 \ar{r} & 0 \ar{r}\ar{d} & \mc{F}_{-1} \ar{r}{id}\ar{d}{\rho_-} & \mc{F}_{-1}\ar{d} \ar{r} & 0 \\
    0 \ar{r} & \mc{F}_0 \ar{d}{id} \ar{r} & \mc{V}_0 \ar{d}{\rho_+} \ar{r} & \mc{F}_0 \ar{r} \ar{d}& 0 \\
    0 \ar{r} & \mc{F}_0 \ar{r} \ar{d} & \mc{F}_1 \ar{r} \ar{d} & \mb{C}_x \ar{r} \ar{d} & 0 \\
      & 0 & 0 & 0,
  \end{tikzcd}
\end{equation}
where the second row induces an extension class $w$ in $\mathrm{Ext}^1(\mc{F}_0,\mc{F}_0)$. We claim that $w$ maps to zero in both $\mathrm{Ext}^1(\mc{F}_0,\mc{F}_1)$ and $\mathrm{Ext}^1(\mc{F}_0,\mc{F}_{-1})$ and that $$d\pi(w,0)=-d\pi(0,w).$$ 
To see this, we consider the following exact sequence of sheaves:
\begin{equation*}
  0\to \mc{F}_{0}/\mc{F}_{-1}\to \mc{F}_1/\mc{F}_{-1}\to \mc{F}_1/\mc{F}_{0}\to 0,
\end{equation*}
which induces a vector $\gamma \in \mathrm{Ext}^1(\mb{C}_x,\mb{C}_x)$. We notice that 
there are two maps from $\mc{F}_{-1}$ to $\mc{V}_0$: the first map is given by $\rho_-$ in the above diagram, and the second map is given by the composition $\mc{F}_{-1}\to \mc{F}_0\to \mc{V}_0$ and is denoted by $\eta_-$. Then we have the following diagram:
\begin{equation*}
  \begin{tikzcd}
     0\ar{r} & \mc{F}_{-1} \ar[r,"0\oplus id"] \ar[d,hook] & \mc{F}_{-1}\oplus \mc{F}_{-1} \ar[r,"id\oplus 0"] \ar[d,"\rho_-\oplus \eta_-",hook] & \mc{F}_{-1} \ar[r] \ar[d,hook]  & 0  \\
     0 \ar{r} & \mc{F}_0 \ar{r}& \mc{V}_0 \ar{r} & \mc{F}_0 \ar{r} & 0.
  \end{tikzcd}
\end{equation*}
Thus, $w$ maps to $0$ in $\mathrm{Ext}^1(\mc{F}_{-1},\mc{F}_0)$, and $\mathrm{coker}(\rho_-\oplus \eta_-
)\cong \mc{F}_{1}/\mc{F}_{-1}$. Hence $d\pi(0,w)=\gamma$. Similarly, we can show that $w$ maps to $0$ in $\mathrm{Ext}^1(\mc{F}_0,\mc{F}_1)$ and $d\pi(w,0)=-\gamma$.

Now we describe the tangent space at a closed point in $\mf{Z}_{-}$, $\mf{Z}_{+}$, and $\mf{Y}$. Given a closed point $(\mc{F}_{0}\subset_{x} \mc{F}_{1}, \mc{F}_{0}'\supset_{y} \mc{F}_{1})$ in $\mf{Z}_{-}$, the tangent space of $\mf{Z}_{-}$ at this point is given by triples 
\begin{equation*}
  \{(w_{1},w_{0},w_{0}')\in \mathrm{Ext}^1(\mc{F}_{1},\mc{F}_{1})\oplus \mathrm{Ext}^1(\mc{F}_{0},\mc{F}_{0})\oplus \mathrm{Ext}^1(\mc{F}_{0}',\mc{F}_{0}')\}
\end{equation*}
such that $w_{1}$ and $w_{0}$ (resp. $w_{1}$ and $w_{0}'$) map to the same element in $\mathrm{Ext}^1(\mc{F}_{1},\mc{F}_{0})$ (resp. $\mathrm{Ext}^1(\mc{F}_{1},\mc{F}_{0}')$) via the maps induced by the inclusion $\mc{F}_{0}\subset \mc{F}_{1}$ (resp. $\mc{F}_{0}'\subset \mc{F}_{1}$). A similar description holds for the tangent space at a closed point in $\mf{Z}_{+}$. By the proof of Proposition 2.39 in \cite{neguct2018hecke}, given a closed point $\{\mc{F}_{-1}\subset_{x}\mc{F}_0\subset_{y} \mc{F}_1, \mc{F}_{-1}\subset_{y}\mc{F}_0'\subset_{x} \mc{F}_1\}\in \mf{Y}$, its tangent space is given by quadruples 
  \begin{equation*}
    (w_{-1},w_{0},w_{0}',w_{1})\in \mathrm{Ext}^1(\mc{F}_{-1},\mc{F}_{-1})\oplus \mathrm{Ext}^1(\mc{F}_{0},\mc{F}_{0})\oplus \mathrm{Ext}^1(\mc{F}_{0}',\mc{F}_{0}')\oplus \mathrm{Ext}^1(\mc{F}_{1},\mc{F}_{1})
  \end{equation*}
  such that $w_{-1}$ and $w_{0}$ (resp. $w_{-1}$ and $w_0'$, $w_0$ and $w_1$, $w_0'$ and $w_1$) map to the same element in $\mathrm{Ext}^1(\mc{F}_{-1},\mc{F}_{0})$ (resp. $\mathrm{Ext}^1(\mc{F}_{-1},\mc{F}_0')$, $\mathrm{Ext}^1(\mc{F}_{0},\mc{F}_1)$ and $\mathrm{Ext}^1(\mc{F}_0',\mc{F}_1)$), with
  $d\pi(w_0,w_1)=d\pi(w_{-1},w_0)$ and $d\pi(w_0',w_1)=d\pi(w_{-1},w_0')$.

\begin{proof}[Proof of \cref{p36}] By \cref{l28} and \cref{eq:1.8}, both $\mc{B}l_{\mf{Z}_{k-1,k}}\mf{Z}_{-}$ and $\mc{B}l_{\mf{Z}_{k,k+1}}\mf{Z}_{+}$ are smooth varieties, and their exceptional divisors are isomorphic to $\mf{Z}_{k-1,k,k+1}^{\bullet}$. By Zariski's main theorem, to show that $\bar{\alpha}_{+}$ and $\bar{\alpha}_{-}$ are isomorphisms, it suffices to show that they are bijections on closed points.

Let $(\mc{F}_{-1}\subset_x\mc{F}_0\subset_x \mc{F}_1)$ be a closed point in $\mf{Z}_{k-1,k,k+1}^{\bullet}$. The tangent space at this point splits as the direct sum of the tangent space of $\mf{Z}_{k-1,k,k+1}^{\bullet}$ and $N$, where $N$ is the normal space generated by the vector $(0,w,-w,0)$, and $w$ is the vector in $\mathrm{Ext}^1(\mc{F}_0,\mc{F}_0)$ defined by \cref{eq:1.10}. 

Now consider a closed point $o:=(\mc{F}_0\subset_x \mc{F}_1)$ in $\mf{Z}_{k-1,k}$. The tangent space $T_o\mf{Z}_-$ splits as $T_o\mf{Z}_{k,k+1}\oplus N_o$, where $N_o$ is the normal space consisting of those $w\in \mathrm{Ext}^1(\mc{F}_0,\mc{F}_0)$ whose image in $\mathrm{Ext}^1(\mc{F}_0,\mc{F}_1)$ is zero. Since $\mathrm{Hom}(\mc{F}_0,\mc{F}_0)\cong \mathrm{Hom}(\mc{F}_0,\mc{F}_1)\cong \mathbb{C}$, we have the left exact sequence
\begin{equation*}
  0\to \mathrm{Hom}(\mc{F}_0,\mb{C}_x)\to \mathrm{Ext}^1(\mc{F}_0,\mc{F}_0)\to \mathrm{Ext}^1(\mc{F}_0,\mc{F}_1).
\end{equation*}
Thus $\mathrm{Hom}(\mc{F}_0,\mb{C}_x)$ is also the normal space of $(\mc{F}_0\subset_x\mc{F}_1)$ in $\mf{Z}_-$. We now prove that $\bar{\alpha}_-$ is a bijection on closed points. For the point $o$, the fiber over this closed point in $\mf{Z}_-$ is parametrized by the one-dimensional vector subspaces of $\mathrm{Hom}(\mc{F}_0,\mathbb{C}_x)$. Given a one-dimensional vector subspace of $\mathrm{Hom}(\mc{F}_0,\mb{C}_x)$ corresponding to the closed point $(\mc{F}_{-1}\subset_x\mc{F}_0\subset_x \mc{F}_1)$ in $\mf{Z}_{k-1,k,k+1}^{\bullet}$, let $w$ be the corresponding vector in $\mathrm{Ext}^1(\mc{F}_0,\mc{F}_0)$. By the above explanation, $w$ is exactly the normal vector of $\mf{Z}_{k-1,k,k+1}^{\bullet}$ in $\mf{Y}$ at the point $(\mc{F}_{-1}\subset_x\mc{F}_0\subset_x \mc{F}_1)$. Hence $\bar{\alpha}_-$ is a bijection on closed points. A similar argument shows that $\bar{\alpha}_+$ is also a bijection on closed points.

\end{proof}

\section{The main theorem}
\label{sec6}
We now state the main theorem of this paper.
\begin{theorem}
  \label{thm61}
  Given integers $i,j,a,b$, there are complexes $\mf{C}_{i,j}^{a,b}\in \D^{b}_{coh}(\mc{M}_{k}\times \mc{M}_{k}\times S\times S)$ and $h_{i}^{a\pm}\in \D^{b}_{coh}(\mc{M}_{k}\times S)$, together with exact triangles
  \begin{align}
    \label{eq:exact1}
  \cdots \to \mf{C}_{i,j}^{a-1,b}\xrightarrow{\mf{c}_{i,j}^{a,b+}}\mf{C}_{i,j}^{a,b}\to R\Delta_{S*}h_{i+j}^{a+}\to \cdots \\
  \label{eq:exact2}
     \cdots \to \mf{C}_{i,j}^{a,b-1}\xrightarrow{\mf{c}_{i,j}^{a,b-}}\mf{C}_{i,j}^{a,b}\to R\Delta_{S*}h_{i+j}^{b-}\to \cdots
  \end{align}
  such that
  \begin{equation*}
      f_ie_j\cong \begin{cases}
        \mf{C}_{i-r,j}^{i+j-r,0} & \text{if } i+j>0, \\
        \mf{C}_{i-r,j}^{i+j-r,i+j} & \text{if } i+j<0,
              \end{cases} \quad 
      R\iota_{*}e_jf_i\cong \begin{cases}
        \mf{C}_{i-r,j}^{i+j-r,0} & \text{if } i+j>0, \\
        \mf{C}_{i-r,j}^{i+j-r,i+j} & \text{if } i+j<0.
              \end{cases}
  \end{equation*}
  Moreover, if $-r<a\leq i+j-r$, then
    \begin{equation}
    \label{eq:neq1}
    h_{i+j-r}^{a+}\cong q^{-a}det(\mc{U}_{k})^{-1}\{\wedge^{r-1+a}\mc{U}_{k}\otimes S^{i+1-a}\mc{U}_{k}\cdots \to S^{i+j}\mc{U}_{k}\}[2-2a-r],
  \end{equation}
  and if $i+j-r<a\leq 0$, then
  \begin{equation}
    \label{eq:neq2}
    h_{i+j-r}^{a-}\cong q^{-a}det(\mc{U}_{k})^{-1}\{\wedge^{-i-j}\mc{U}_{k}^{\vee}\to \cdots \to \wedge^{-i-j+a}\mc{U}_{k}^{\vee}\otimes S^{-a}\mc{U}_{k}^{\vee}\}[1-a-r]. 
  \end{equation}
  When $i+j=0$, we have the decomposition
    $$f_{i}e_{j}\cong R\iota_{*}e_{j}f_{i}\bigoplus_{a=-r+1}^{0} R\Delta_{S*}(q^{-a}det(\mc{U}_{k})^{-1}\mc{O}_{\mc{M}_{k}\times S})[1-2a-r].$$
\end{theorem}

\subsection{The complexes $h_{i}^{a\pm}$ and $\mf{C}_{i,j}^{a,b}$}
We recall the resolution of the universal sheaves in \cref{eq:1.7}. By \cref{p25}, the moduli space $\mf{Z}_+\subset \mb{P}_{\mf{Z}_{k,k+1}\times S}(q\mc{W}_{k+1}^{\vee})$ 
is the zero locus of the cosection $\taut_{q\psi_{k+1}^{\vee}}$. Let $\Bl_+:=\Bl_{\mf{Z}_{k,k+1}}\mb{P}_{\mf{Z}_{k,k+1}\times S}(q\mc{W}_{k+1}^{\vee})$, with exceptional divisor $D:=\mb{P}_{\mf{Z}_{k,k+1}}(q\mc{L}\mc{W}_{k}^{\vee})$ and projection morphism
$$pr_{+}:\Bl_{+}\to \mb{P}_{\mf{Z}_{k,k+1}\times S}(q\mc{W}_{k+1}^{\vee}).$$
Let $\tau_D:D\to \mb{P}_{\mf{Z}_{k,k+1}\times S}(q\mc{W}_{k+1}^{\vee})$ be the closed embedding of the exceptional divisor. We define the line bundle $\mc{L}_1$ on $\Bl_+$ as the pullback of $\mc{L}$ from $\mf{Z}_{k,k+1}$, and we define the line bundle $\mc{L}_1'$ on $\Bl_+$ as the pullback of $\mc{O}_{\mb{P}_{\mf{Z}_{k,k+1}\times S}(q\mc{W}_{k+1}^{\vee})}(-1)$. Let
\begin{equation}
  \label{eq31}
  \taut_+:q\mc{L}_{1}'\mc{V}_{k+1}^{\vee}\to \mc{O}_{\Bl_{+}}  
\end{equation}
be the pullback of $\taut_{q\psi_{k+1}^{\vee}}$ to $\Bl_+$, and let $\mb{R}\overline{\mf{Z}}_{+}$ be the derived zero locus of $\taut_+$. By \cref{l27}, $\taut_+$ induces a morphism of locally free sheaves on $\Bl_+$,
\begin{equation}
  \label{eq32}
  \widetilde{\taut}_+:q\mc{L}_{1}'\mc{V}_{k+1}^{\vee}\otimes \mc{O}(D)\to \mc{O}_{\Bl_{+}},
\end{equation}
whose zero locus is $\mf{Y}$ by \cref{p36}. Moreover, the definitions of the line bundles $\mc{L}_1$ and $\mc{L}_1'$ on $\Bl_+$ are compatible with the corresponding line bundles on $\mf{Y}$ under pullback. Given integers $i,j$ and $a\leq b$, we define the complex of coherent sheaves $\mf{B}_{i,j}^{a,b}$ on $\Bl_{+}$ whose degree-$m$ term is
  \begin{equation}
    \label{d64}
     (\mf{B}^{a,b})_{m}=
    \begin{cases}
      \wedge^{-m}(q\mc{L}_{1}' \mc{V}_{k+1}^{\vee})\otimes \mc{O}((a-m)D) & m<a, \\
      \wedge^{-m}(q\mc{L}_{1}' \mc{V}_{k+1}^{\vee}) & a\leq m\leq b, \\
     \wedge^{-m}(q\mc{L}_{1}' \mc{V}_{k+1}^{\vee})\otimes \mc{O}((b-m)D) & m>b.
    \end{cases}
  \end{equation}
  The differential $d_{m}:(\mf{B}^{a,b})_{m}\to (\mf{B}^{a,b})_{m+1}$ is induced by $\wedge^{-m}\taut_+$ when $a\leq m<b$, and by $\wedge^{-m}\widetilde{\taut}_+$ when $m<a$ or $m\geq b$; these maps are defined in \cref{eq31} and \cref{eq32}. For any integer $a$, we have $\mf{B}^{a,0}\cong \mf{B}^{a,1}\cong \mf{B}^{a,2}=\cdots$. For any integer $b$, we have $\mf{B}^{a,b}=\mf{B}^{a-1,b}$ when $a$ is sufficiently small; we denote this complex by $\mf{B}^{-\infty,b}$. Given any pairs of integers $(a_{1}\leq b_{1})$ and $(a_{2}\leq b_{2})$, there is a canonical morphism
  $$\mf{b}^{(a_{1},b_{1}),(a_{2},b_{2})}:\mf{B}^{a_{1},b_{1}}\to \mf{B}^{a_{2},b_{2}}$$ 
  induced by inclusion. By \cref{m57}, $\mf{B}^{a,b}$ and $\mf{b}^{(a_{1},b_{1}),(a_{2},b_{2})}$ are respectively objects and morphisms in $\D^{b}_{coh}(\mb{R}\bar{\mf{Z}}_{+})$.

Let $\mc{O}_{D}(1):=\mc{O}_{\mb{P}_{\mf{Z}_{k,k+1}}(q\mc{L}\mc{W}_{k}^{\vee})}(1)$ be the tautological line bundle on $D$. We truncate the Koszul complex of $\taut_{q\mc{L}\psi_{k,k+1}^{\vee}}$ to define complexes $\mf{H}^{a+}$ and $\mf{H}^{a-}$ of coherent sheaves on $D$ whose degree-$m$ terms are
  \begin{align}
    \label{eq:33}
     (\mf{H}^{a+})_{m}=
    \begin{cases}
      0 & m\geq a, \\
      \wedge^{-m}(q\mc{L}\mc{V}_{k+1}^{\vee})\mc{O}_{D}(m-a) & m<a;
    \end{cases}\\
      \label{eq:34}
    (\mf{H}^{a-})_{m}=
    \begin{cases}
      \wedge^{-m}(q\mc{L}\mc{V}_{k+1}^{\vee})\mc{O}_{D}(m-a) & m\geq a, \\
      0 & m< a.
    \end{cases}
  \end{align}
We define $\mf{H}_{i}^{a\pm}:=\mc{L}^{i}\mf{H}^{a\pm}$ and $h_{i}^{a\pm}:=R((p_{k}\times \pi_{k})\circ pr_{q\mc{L}\mc{W}_{k}^{\vee}})_{*}(\mf{H}_{i}^{a\pm})$. By the definitions in \cref{d64,eq:33,eq:34}, we have
\begin{equation}
  \label{eq:exact3}
\mf{B}_{i,j}^{a,b}/\mf{B}_{i,j}^{a-1,b}=\tau_{D*}\mf{H}_{i+j}^{a+} \quad \text{and} \quad \mf{B}_{i,j}^{a,b}/\mf{B}_{i,j}^{a,b-1}=\tau_{D*}\mf{H}_{i+j}^{b-}.
\end{equation}

Let $\Delta_{k+1}$ be the diagonal embedding of $\mc{M}_{k+1}$, and let
$$\tau_{k}:\mf{Z}_{k,k+1}\to \mb{P}_{\mc{M}_{k+1}\times S}(q\mc{W}_{k+1}^{\vee})$$
be the closed embedding obtained by composing the closed embedding $\Delta_{+}$ with the closed embedding of $\mf{Z}_+$ into $\mb{P}_{\mc{M}_{k+1}\times S}(q\mc{W}_{k+1}^{\vee})$. We have the following diagram, in which all squares are derived Cartesian diagrams:
\begin{equation}
  \label{cref418}
  \begin{tikzcd}
 \mb{R}\bar{\mf{Z}}_{+} \ar{d} \ar{rr}{\mc{z}}&       & \mf{Z}_{k,k+1}\times \mf{Z}_{k,k+1} \ar{d}{id\times \tau_{k}} \\
\Bl_{+} \ar{r}{pr_{+}} & \mb{P}_{\mf{Z}_{k,k+1}\times S}(q\mc{W}_{k+1}^{\vee}) \ar{r}{\gamma}\ar{d} & \mf{Z}_{k,k+1}\times \mb{P}_{\mc{M}_{k+1}\times S}(q\mc{W}_{k+1}^{\vee})\ar{d}{q_k\times pr_{q\mc{L}\mc{W}_{k+1}^{\vee}}} \\
& \mc{M}_{k+1}\ar{r}{\Delta_{k+1}} & \mc{M}_{k+1}\times \mc{M}_{k+1}
  \end{tikzcd}
\end{equation}
where $\gamma$ is the pullback of $\Delta_{k+1}$ and $\mc{z}$ is the pullback of $\gamma \circ pr_{+}$. Consider the morphism
\begin{align*}
  \pi:=p_{k}^{1}\times p_{k}^{2}\times \pi_{k}^{1}\times \pi_{k}^{2}:\mf{Z}_{k,k+1}\times \mf{Z}_{k,k+1}&\to \mc{M}_{k}\times \mc{M}_{k}\times S\times S,
\end{align*}
which maps pairs $(\mc{F}_{-1}\subset_{x} \mc{F}_{0}, \mc{F}_{-1}'\subset_{y} \mc{F}_{0}')$ to $(\mc{F}_{0},\mc{F}_{0}',x,y)$. Define
$$\mf{C}_{i,j}^{a,b}:=R(\pi\circ \mc{z})_{*}(\mc{L}_1^i\mc{L}_1'^j\otimes \mf{B}^{a,b}), \quad \mf{c}_{i,j}^{(a_{1},b_{1}),(a_{2},b_{2})}:=R(\pi\circ \mc{z})_{*}(\mc{L}_1^i\mc{L}_1'^j\otimes \mf{b}^{(a_{1},b_{1}),(a_{2},b_{2})})$$
for any integers $a,b,i,j$, and set $\mf{c}_{i,j}^{a,b+}:=\mf{c}_{i,j}^{(a-1,b),(a,b)}$ and $\mf{c}_{i,j}^{a,b-}:=\mf{c}_{i,j}^{(a,b-1),(a,b)}$.
Let $\tau_{D}$ be the closed embedding from $D$ to $\Bl_{+}$. Recalling $\mf{H}_{i}^{a\pm}$ from \cref{eq:33}, we have
$$\pi\circ \mc{z}\circ \tau_{D}=\Delta_{S} \circ (p_{k}\times \pi_{k})\circ pr_{q\mc{L}\otimes \mc{W}_{k}^{\vee}},$$
and hence $R(\pi\circ \mc{z}\circ \tau_{D})_{*}\mf{H}_{i}^{a\pm}=R\Delta_{S*}h_{i}^{a\pm}.$ Therefore, \cref{eq:exact3} gives the exact triangles \cref{eq:exact1} and \cref{eq:exact2}.

\subsection{The explicit formula for $h_{i}^{a\pm}$ and the proof of \cref{thm61}}
\label{sec4} By \cref{lemma:n39}, the truncations in \cref{eq:33,eq:34} induce a canonical morphism $\mb{H}_{i}^{a}:\mf{H}_{i}^{a+}[-1]\to \mf{H}_{i}^{a-}$
such that its cone resolves the sheaf $\mc{L}_{2}^{a}\mc{L}^{i-a}_{1}\mc{O}_{\Delta_{\mf{Y}}}$. By \cref{n510}, the morphism
\begin{equation}
  \label{n55}
  (R((p_{k}\times \pi_{k})\circ pr_{q\mc{L}\otimes \mc{W}_{k}^{\vee}})_{*}\mb{H}_{i}^{a}):h_{i}^{a+}[-1]\to h_{i}^{a-}  
\end{equation}
is an isomorphism when $-r\leq i-a<0$. We now compute $h_{i}^{a\pm}$.
\begin{proposition}
  \label{lcor} \label{n54} \label{n56} Equations \cref{eq:neq1} and \cref{eq:neq2} hold. Moreover, $h_{i}^{a+}\cong0$ if $a\leq -r$, and $h_{i}^{a-}\cong0$ if $i<a\leq i+r$ and $i+r\neq 0$. When $i=-r$ and $-r<a\leq 0$, we have
  \begin{equation}
    \label{eq:neq3}
   h_{i}^{a-}\cong q^{-a}det(\mc{U}_{k})^{-1}\mc{O}_{\mc{M}_{k}\times S}[1-2a-r].
 \end{equation}
 Moreover, \cref{eq:neq3} is the image of the truncation morphism $\mf{h}_{-r}^{a-}\to \wedge^{-a}(q\mc{L}\mc{V}_{k+1})^{\vee}[-a]$ under the functor $R((p_{k}\times \pi_{k})\circ pr_{q\mc{L}\otimes \mc{W}_{k}^{\vee}})_{*}$:
 \begin{equation}
   \label{eq:4115}
    h_{i}^{a-}\cong R((p_{k}\times \pi_{k})\circ pr_{q\mc{L}\otimes \mc{W}_{k}^{\vee}})_{*}(\mc{L}^{-r}\wedge^{-a}(q\mc{L}\mc{V}_{k+1})^{\vee}[-a]).
 \end{equation}
\end{proposition}
\begin{proof}
The short exact sequence
  $$0\to \mc{L}^{-1}\mc{V}_{k+1}\to \mc{L}^{-1}\mc{V}_{k}\to \mc{O}_{D}\to 0$$
  on $D$ induces the following long exact sequences on $D$ for any $m>0$:
\begin{align}
  \label{lequa4.5}  0 \to \mc{O} \to \mc{L}\mc{V}_{k}^{\vee} \to \cdots \to \wedge^{m}(\mc{L}\mc{V}_{k}^{\vee})\to \wedge^{m}(\mc{L}\mc{V}_{k+1}^{\vee})\to 0,\\
 \label{lequa4.6} 0\to \wedge^{m}(\mc{L}^{-1}\mc{V}_{k+1})\to \wedge^{m}(\mc{L}^{-1}\mc{V}_{k})\to \cdots \to \mc{L}^{-1}\mc{V}_{k}\to \mc{O}\to 0. 
\end{align}
Recalling the morphism $\mc{L}^{-1}\psi_{k,k+1}:\mc{W}_{k}\to \mc{V}_{k+1}$ from \cref{psik}, \cref{lequa4.5} and \cref{lequa4.6} give
\begin{align}
    \label{eq311}
  S^{m}(\mc{L}\otimes \psi_{k,k+1}^{\vee})&\cong \{\mc{O}\to \cdots \to S^{m-1}(\mc{L}\mc{U}_{k}^{\vee})\to S^{m}(\mc{L}\mc{U}_{k}^{\vee})\}[m],\\
  \label{eq312}
  \wedge^{m}(\mc{L}^{-1}\otimes \psi_{k,k+1})&\cong \{\wedge^{m}(\mc{L}^{-1}\mc{U}_{k})\to \wedge^{m-1}(\mc{L}^{-1}\mc{U}_{k})\to \cdots \to \mc{O}\}[-m].
\end{align}

By \cref{hart}, if $a\leq -r$, then $Rpr_{w_{k}*}\mf{H}^{a+}\cong0$, and hence $h_{i}^{a+}\cong0$. When $a> -r$, \cref{eq312} gives
\begin{align*}
  Rpr_{w_{k}*}\mf{H}^{a+}
                         &\cong det(q\mc{L}\mc{V}_{k+1}^{\vee})det(q\mc{L}\mc{W}_{k}^{\vee})^{-1}\wedge^{rank(\mc{V}_{k+1})-rank(\mc{W}_{k})+a}(q^{-1}\mc{L}^{-1}\psi_{k,k+1})[1-a] \\
                         &\cong q^{r-1}det(\mc{L}\mc{V}_{k+1}^{\vee})det(\mc{L}\mc{W}_{k}^{\vee})^{-1}q^{-r+1-a}\wedge^{r-1+a}(\mc{L}^{-1}\psi_{k,k+1})[1-a] \\
                         &\cong q^{-a}\mc{L}^{r}det(\mc{U}_{k})^{-1}\wedge^{r-1+a}(\mc{L}^{-1}\psi_{k,k+1})[1-a] \\
                       &\cong q^{-a}\mc{L}^{r}det(\mc{U}_{k})^{-1}\{\wedge^{r-1+a}(\mc{L}^{-1}\mc{U}_{k})\to \cdots \to \mc{L}^{-1}\mc{U}_{k}\to \mc{O}\}[2-2a-r]. 
\end{align*}
Thus \cref{eq:neq1} follows from \cref{cor:5.3ext}. When $a\leq 0$, by \cref{hart} and \cref{eq311}, we have
\begin{equation}
  \label{714}
  Rpr_{w_{k}*}\mf{H}^{a-}\cong q^{-a}S^{-a}(\mc{L}\psi_{k,k+1}^{\vee})\cong q^{-a}\{\mc{O}\to \cdots \to S^{-a}(\mc{L}\mc{U}_{k}^{\vee})\}[-a].
\end{equation}
This gives \cref{eq:neq2} and \cref{eq:neq3}. When $i<a\leq i+r$ and $i+r\neq 0$, \cref{714} gives
$$h_{i}^{a-}\cong q^{-a}det(\mc{U}_{k})^{-1}\{\wedge^{i+r}(\mc{U}_{k}^{\vee})\to \cdots \to S^{i+r}(\mc{U}_{k}^{\vee})\}\cong0.$$

It remains to prove \cref{eq:4115}. By \cref{lequa4.5},
  \begin{align*}
    Rpr_{q\mc{L}\mc{W}_{k}^{\vee}}(\mc{L}^{-r}\wedge^{-a}(q\mc{L}\mc{V}_{k+1})^{\vee}[-a])&\cong \mc{L}^{-r}\wedge^{-a}(q\mc{L}\mc{V}_{k+1})^{\vee}[-a] \\
    &\cong q^{-a}\{\mc{O}\to \mc{L}\mc{V}_{k}^{\vee}\to \cdots \to \wedge^{-a}(\mc{L}\mc{V}_{k}^{\vee})\}[-a].
  \end{align*}
The truncation morphism $S^{m}(\mc{L}\mc{U}_{k}^{\vee})\to \wedge^{m}(\mc{L}\mc{V}_{k}^{\vee})$
  induces the morphism
  $$R(pr_{q\mc{L}\mc{W}_{k}^{\vee}})_{*}\mf{h}_{-r}^{a-}\to R(pr_{q\mc{L}\mc{W}_{k}^{\vee}})_{*}(\mc{L}^{-r}\wedge^{-a}(q\mc{L}\mc{V}_{k+1})^{\vee}[-a]).$$
Hence, by \cref{cor:5.3ext} and \cref{n56}, we obtain the isomorphism \cref{eq:4115}.
\end{proof}

\begin{proof}[Proof of \cref{thm61}]
  By \cref{d64} and \cref{p36}, as objects of $D_{coh}^{b}(\mb{R}\bar{\mf{Z}}_{+})$, we have $\mc{L}_{1}^{i}\mc{L}_{1}'^{j}\mc{O}(lD)\otimes \mc{O}_{\mf{Y}}\cong \mf{B}_{i,j}^{l,l}$. Thus, by \cref{n42},
\begin{equation}
  \label{eq:6022}
  R\iota_{*}e_{j}f_{i}\cong R(\pi\circ \mc{z})_{*}(\mf{B}_{i-r,j}^{i+j-r,i+j-r})\cong \mf{C}_{i-r,j}^{i+j-r,i+j-r}.
\end{equation}
On the other hand, we have
\begin{align}
\label{eq:6021}  f_{i}e_{j}&\cong R\pi_{*}\circ L(id\times \tau_{k})^{*}\circ R\gamma_{*}(\mc{L}^{i-r}_{1}\mc{L}_{1}'^{j}) & \\
            &\cong R\pi_{*}L(id\times \tau_{k})^{*}R(\gamma\circ pr_{\mf{t}_{k}})_{*}(\mc{L}^{i-r}_{1}\mc{L}_{1}'^{j}) \nonumber \\
            &\cong R(\pi\circ \mc{z})_{*}(\mf{B}_{i-r,j}^{-\infty,0})\cong \mf{C}_{i-r,j}^{-\infty,0}. & \text{by \cref{n613}} \nonumber
\end{align}
By \cref{n56}, $h_{i+j-r}^{a+}\cong 0$ when $a\leq -r$. Thus, for any $l\leq -r$, we have $\mf{C}_{i-r,j}^{-\infty,0}\cong \mf{C}_{i-r,j}^{l,0}$; in particular,
\begin{equation*}
  \mf{C}_{i-r,j}^{-\infty,0}\cong \begin{cases}
    \mf{C}_{i-r,j}^{-r,0} & \text{if } i+j>0, \\
    \mf{C}_{i-r,j}^{i+j-r,0} & \text{if } i+j<0.
  \end{cases}
\end{equation*}
On the other hand, if $i+j\neq 0$, then again by \cref{n56}, $h_{i+j-r}^{a-}\cong 0$ when $i+j-r<a\leq i+j$. Thus $\mf{C}_{i-r,j}^{i+j-r,i+j-r}\cong \mf{C}_{i-r,j}^{i+j-r,l}$ for all $i+j-r\leq l\leq i+j$ if $i+j\neq 0$; in particular,
\begin{equation*}
  \mf{C}_{i-r,j}^{i+j-r,i+j-r}\cong \begin{cases}
    \mf{C}_{i-r,j}^{i+j-r,0} & \text{if } i+j>0, \\
    \mf{C}_{i-r,j}^{i+j-r,i+j} & \text{if } i+j<0.
  \end{cases}
\end{equation*}

The preceding arguments prove the case $i+j\neq 0$. When $i+j=0$, we have the exact triangle
 \begin{equation}
   \label{eq:exacttriangle}
   \cdots \to \mf{C}_{i-r,j}^{-r,-r} \xrightarrow{\mf{c}_{i-r,j}^{(-r,-r),(-r,0)}} \mf{C}_{i-r,j}^{-r,0}\to R(\pi\circ \mc{z})_{*}(\mf{B}_{i-r,j}^{-r,0}/\mf{B}_{i-r,j}^{-r,-r})\to \cdots.
 \end{equation}
By \cref{n55}, we have 
$R\pi_{*}R\mc{z}_{*}(\mf{B}_{i-r,j}^{0,0}/\mf{B}_{i-r,j}^{-r,-r})=0$, and hence
\begin{equation}
  \label{laeq}
\mf{c}_{i-r,j}^{(-r,-r),(0,0)}=\mf{c}_{i-r,j}^{(-r,0),(0,0)}\circ \mf{c}_{i-r,j}^{(-r,-r),(-r,0)}:\mf{C}_{i-r,j}^{-r,-r}\to \mf{C}_{i-r,j}^{0,0}  
\end{equation}
is an isomorphism. Hence $\mf{c}_{i-r,j}^{(-r,-r),(0,0)}$ has a left inverse and \cref{eq:exacttriangle} splits. Moreover, the restriction morphism $\mc{L}_{1}^{i-r}\mc{L}_{1}'^{j}\wedge^{-a}(q\mc{L}_{1}'\mc{V}_{k+1}^{\vee})\to \mc{L}^{-r}\wedge^{-a}(q\mc{L}_{1}'\mc{V}_{k+1}^{\vee})\otimes \mc{O}_{D}$ induces a morphism of dg-modules $\mf{B}_{i-r,j}^{-r,0}\to \oplus_{a=-r}^{0}\mc{L}^{-r}\wedge^{-a}(q\mc{L}_{1}'\mc{V}_{k+1}^{\vee})\otimes \mc{O}_{D}$
such that the image of $\mf{B}_{i-r,j}^{-r,-r}$ is $0$. Hence it induces a morphism
$$\mf{B}_{i-r,j}^{-r,0}/\mf{B}_{i-r,j}^{-r,-r}\to \bigoplus_{a=-r}^{0}\mc{L}^{-r}\wedge^{-a}(q\mc{L}_{1}'\mc{V}_{k+1}^{\vee})\otimes \mc{O}_{D}.$$
By \cref{lcor}, we have
\begin{align*}
  R(\pi\circ \mc{z})_{*}(\mf{B}_{i-r,j}^{-r,0}/\mf{B}_{i-r,j}^{-r,-r})&\cong \bigoplus_{a=-r}^{0}R(\pi\circ \mc{z})_{*}\mc{L}^{-r}\wedge^{-a}(q\mc{L}_{1}'\mc{V}_{k+1}^{\vee})\otimes \mc{O}_{D}. \\
  &\cong \bigoplus_{a=-r}^{0} R\Delta_{S*}(q^{-a}det(\mc{U}_{k})^{-1}\mc{O}_{\mc{M}_{k}\times S})[1-2a-r].
\end{align*}
\end{proof}

\end{document}